\documentclass[10pt,a4paper]{amsart}
\usepackage{amssymb}
\numberwithin{equation}{section}
\newtheorem{theorem}{Theorem}[section]

\newtheorem{proposition}[theorem]{Proposition}
\newtheorem{lemma}[theorem]{Lemma}
\theoremstyle{remark}
\newtheorem{remark}{Remark}[section]

\theoremstyle{definition}

\newcommand{\triple}[1]{{|\!|\!|#1|\!|\!|}}
\newcommand{\R}{\mathbb{R}}
\newcommand{\C}{\mathbb{C}}
\newcommand{\D}{\mathcal{D}}
\newcommand{\xx}{\langle x\rangle}
\newcommand{\yy}{\langle y\rangle}
\newcommand{\zz}{\langle z\rangle}
\newcommand{\DD}{\langle D\rangle}
\newcommand{\DDW}{\langle D_{W}\rangle}

\begin{document}

\title%
[Strichartz and smoothing estimates]%
{Strichartz and smoothing
estimates for dispersive equations with magnetic potentials}
\begin{abstract}
  We prove global smoothing and Strichartz estimates for the
  Schr\"o\-din\-ger, wave, Klein-Gordon equations and for
  the massless and massive Dirac systems, perturbed with singular
  electromagnetic potentials. We impose a smallness condition
  on the magnetic part, while the electric part can be large.
  The decay and regularity assumptions on the
  coefficients are close to critical.
\end{abstract}
\date{\today}    
\author{Piero D'Ancona}
\address{Piero D'Ancona:
Unversit\`a di Roma ``La Sapienza'',
Dipartimento di Matematica,
Piazzale A.~Moro 2, I-00185 Roma, Italy}
\email{dancona@mat.uniroma1.it}

\author{Luca Fanelli}
\address{Luca Fanelli:
Unversit\`a di Roma ``La Sapienza'',
Dipartimento di Matematica,
Piazzale A.~Moro 2, I-00185 Roma, Italy}
\email{fanelli@mat.uniroma1.it}

\maketitle

\section{Introduction}\label{sec.introd}

Strichartz estimates have become a standard
tool in the study of linear and nonlinear evolution equations.
They are available for a large class of constant
coefficient equations, by the methods of
\cite{GinibreVelo95-generstric} and \cite{KeelTao98-endpoinstric}.
In a sense, they represent the modern energy estimates,
and are especially effective for problems of low regularity
and global existence for nonlinear equations.

Using the notations 
$L^{p}L^{q}=L^{p}(\mathbb{R}_{t};L^{q}(\mathbb{R}^{n}_{x}))$,
$\|f\|\lesssim\|g\|$ to mean $\|f\|\le C\|g\|$,
and $H_{q}^{s}$ and $\dot H_{q}^{s}$ to denote the spaces
with norms
\begin{equation*}
  \|f\|_{\dot H_{q}^{s}}=\|\DD^{s}f\|_{L^{q}},\qquad
  \|f\|_{\dot H_{q}^{s}}=\||D|^{s}f\|_{L^{q}}.
\end{equation*}
where $\DD=(1-\Delta)^{1/2}$, $|D|=(-\Delta)^{1/2}$,
the Strichartz estimates for the Schr\"odinger equation 
take the following form: for $n\ge2$,
\begin{equation*}
  \|e^{it \Delta}f\|_{L^{p}L^{q}}\lesssim\|f\|_{L^{2}},
\end{equation*}
provided the couple $(p,q)$ is \emph{Schr\"odinger admissible}:
\begin{equation}\label{eq.schroadmi}
  \frac2p+\frac nq=\frac n2,\qquad
  2\le p\le \infty,\qquad
  \frac{2n}{n-2}\ge q\ge 2,\qquad
  q\neq\infty.
\end{equation}
The couple $(p,q)=(2,2n/n-2)$ is called the \emph{endpoint}
and is allowed when $n>2$.

For the wave equation the estimates can be written as follows:
for $n\ge3$,
\begin{equation*}
  \|e^{it|D|}f\|
  _{L^{p}\dot H_{q}^{\frac 1q-\frac1p-\frac12}}
  \lesssim\|f\|_{L^{2}},
\end{equation*}
provided the couple $(p,q)$ is \emph{wave admissible}:
\begin{equation}\label{eq.waveadmis}
  \frac2p+\frac{n-1}{q}=\frac{n-1}{2},\qquad
  2\le p\le \infty,\qquad
  \frac{2(n-1)}{n-3}\ge q\ge 2,\qquad
  q\neq\infty.
\end{equation}
The wave equation endpoint is $(p,q)=(2,2(n-1)/(n-3))$
and is allowed in dimension $n>3$.

Finally for the Klein-Gordon equation we have: for $n\ge2$,
\begin{equation*}
  \|e^{it\DD}f\|
  _{L^{p} H_{q}^{\frac 1q-\frac1p-\frac12}}
  \lesssim\|f\|_{L^{2}},
\end{equation*}
provided $(p,q)$ is Schr\"odinger admissible
(see the Appendix for a proof of the last estimate,
for which a reference is not immediately available).

We shall also be interested in the decay properties of
the Dirac equation, which is a $4\times4$ constant
coefficient system of the form
\begin{equation*}
  iu_{t}+\D u=0
\end{equation*}
in the \emph{massless} case, and
\begin{equation*}
  iu_{t}+\D u+{\beta}u=0
\end{equation*}
in the \emph{massive} case. Here
$u:\mathbb{R}_{t}\times \mathbb{R}^{3}_{x}\to \mathbb{C}^{4}$,
the operator $\D$ is defined as
\begin{equation*}
  \D=\frac1i\sum_{k=1}^{3}\alpha_{k}\partial_{k}
\end{equation*}
and the $4\times4$ \emph{Dirac matrices} can be written
\begin{equation*}
  \alpha_{k}=
  \begin{pmatrix}
  0 &\sigma_{k}  \\
  \sigma_{k} & 0
  \end{pmatrix},\quad
  \beta=
  \begin{pmatrix}
  I_{2} & 0\\
   0& -I_{2}
  \end{pmatrix},\qquad
  k=1,2,3
\end{equation*}
in terms of the \emph{Pauli matrices}
\begin{equation*}
  I_{2}=
  \begin{pmatrix}
  1 &0 \\
  0 &1
  \end{pmatrix},\quad
  \sigma_{1}=
  \begin{pmatrix}
  0 &1 \\
  1 &0
  \end{pmatrix},\quad
  \sigma_{2}=
  \begin{pmatrix}
  0 &-i \\
  i &0
  \end{pmatrix},\quad
  \sigma_{3}=
  \begin{pmatrix}
  1 &0 \\
  0 &-1
  \end{pmatrix}.
\end{equation*}
Then the solution $u(t,x)=e^{it\D}f$ of the massless Dirac system
with initial value $u(0,x)=f(x)$ satisfies the Strichartz
estimate:
\begin{equation*}
  \|e^{it\D}f\|
  _{L^{p}\dot H_{q}^{\frac 1q-\frac1p-\frac12}}
  \lesssim\|f\|_{L^{2}},\qquad
  n=3,
\end{equation*}
for all wave admissible $(p,q)$,
while in the massive case we have
\begin{equation*}
  \|e^{it(\D+\beta)}f\|
  _{L^{p} H_{q}^{\frac 1q-\frac1p-\frac12}}
  \lesssim\|f\|_{L^{2}},\qquad
  n=3,
\end{equation*}
for all Schr\"odinger admissible $(p,q)$
(see the Appendix for more details).

In view of the applications, it is an important problem
to extend Strichartz estimates to more general
equations with
variable coefficients, possibly of low regularity
in order to retain
the advantages over classical energy methods. Indeed,
in recent years a large number of works
have investigated this kind of problem.
In the case of potential perturbations like
\begin{equation*}
  iu_t-\Delta u+V(x)u=0, \qquad
  \square u+V(x)u=0,
\end{equation*}
Strichartz estimates are now
fairly well understood.
We mention among the many works 
\cite{BurqPlanchonStalker04-a},
\cite{GoldbergSchlag04-a},
\cite{Goldberg04-a},
\cite{RodnianskiSchlag04-a},
\cite{Schlag05-b}
and the survey
\cite{Schlag05-a}
for the Schr\"odinger equation, and
\cite{Cuccagna00-b},
\cite{GeorgievVisciglia03-a},
\cite{DAnconaPierfelice05-a}
for the wave equation.
We also mention the wave operator approach of Yajima
(\cite{Yajima95-b}, \cite{Yajima95-a},
\cite{Yajima99-a}, \cite{ArtbazarYajima00}),
which was recently optimized
in dimension 1 in \cite{DAnconaFanelli06-b}.

Results are much less complete
in the case of first order
perturbations i.e. \emph{magnetic potentials}
\begin{equation*}
  iu_t+\Delta u+a\cdot \nabla u+bu=0,\qquad
  \square u+a\cdot \nabla u+bu=0.
\end{equation*}
Concerning Strichartz estimates
for the Schr\"odinger equation with small potentials
$a,b$ we recall at least the papers
\cite{stef}, \cite{GST}; in 3D the recent
work \cite{ErdoganGoldbergSchlag-a} handles for the
first time the case of large magnetic potentials.
For the wave equation with small magnetic potentials,
partial Strichartz estimates
were obtained in 3D in \cite{CuccagnaSchirmer01-a}
in the case of smooth, rapidly decaying coefficients.
The dispersive estimate in 3D was proved in
\cite{DAnconaFanelli06-c} for the magnetic wave equation
with small singular potentials and for the massless
Dirac system with a small singular matrix potential.
We must also mention the papers
\cite{StaffilaniTataru02-a},
\cite{RobbianoZuily-a},
\cite{TataruGeba05-a}
containig some local estimates in
the fully variable coefficient case. Only in the
one dimensional case the optimal dispersive estimates
for the case of fully variable singular coefficients
have been proved in \cite{DAnconaFanelli06-b}.

A method of proof which is very efficient in
the case of electric potentials
was introduced in \cite{RodnianskiSchlag04-a}
and further developed in \cite{BurqPlanchonStalker04-a}.
The main idea is to combine Strichartz estimates for the
free equation with Kato smoothing estimates for
the perturbed equation. The
same method is used in
\cite{ErdoganGoldbergSchlag-a} for the
3D Schr\"odinger equation with a large magnetic potential.

Our goal here is to apply a suitable modification
of this method in a systematic way to several
equations perturbed with magnetic potentials:
Schr\"odinger, wave and Klein-Gordon equations,
and the Dirac system with and without mass.

Thus consider a magnetic Schr\"odinger operator
\begin{equation}\label{eq:magn1}
  H=-(\nabla+iA(x))^{2}+B(x),
\end{equation}
which is selfadjoint under the following assumptions:
$A_{j}$ and $B$ are real valued, and
\begin{equation}\label{eq.assco2}
    \|B\|_{L^{n/2,\infty}}<\infty,\qquad
    \|B_{-}\|_{L^{n/2,\infty}}<\delta,\qquad
    \|A\|_{L^{n,\infty}}<\delta
\end{equation}
for some $\delta$ sufficiently small
(see Lemma \ref{lem.selfadj} below).
Here $L^{p,\infty}=L^{p}_{w}$ denotes the
Lorentz or weak Lebesgue
space. However, in order to state our results, it is more
convenient to represent the operator in the form
\begin{equation}\label{eq:magn2}
  H \equiv-\Delta+W(x,D)\equiv-\Delta+a(x)\cdot \nabla +b(x)
\end{equation}
and to make the abstract assumption that $H$
is selfadjoint. In view of \eqref{eq.assco2},
the following explicit conditions on $a,b$
are sufficient (but not necessary)
for the selfadjointness of $H$:
\begin{equation}\label{eq:self}
  a(x)\text{ is pure imaginary,}\quad
  \Im b=-i \nabla \cdot a
\end{equation}
and
\begin{equation}\label{eq.assco3}
    \|\nabla a\|_{L^{n/2,\infty}}+\|b\|_{L^{n/2,\infty}}<\infty,\qquad
    \|\Re b_{-}\|_{L^{n/2,\infty}}<\delta,\qquad
    \|a\|_{L^{n,\infty}}<\delta
\end{equation}
for a small enough $\delta$.

Our first result concerns smoothing estimates of Kato-Yajima
type for the scalar Schr\"odinger, wave and
Klein-Gordon equations. Besides being a necessary tool to prove
the Strichartz estimates, they have also an independent
interest (see e.g. \cite{benar}, \cite{Kato66}, \cite{KY}).
Notice in particular that we allow a
singularity at 0 in the coefficient, and that
the electric potential can be large, while the
magnetic term must satisfy a smallness condition. We shall use the
following weight functions:
\begin{equation*}
  \tau_\epsilon(x)=
  \begin{cases}
      |x|^{\frac12-\epsilon}+|x|&\text{if $n\ge3$,}\\
      |x|^{\frac12-\epsilon}+|x|^{1+\epsilon}&\text{if $n=2$}
  \end{cases}
\end{equation*}
and
\begin{equation*}
  w_\sigma(x)=|x|(1+|\log|x||)^\sigma,\qquad
  \sigma>1.
\end{equation*}
Then we have:

\begin{proposition}[Smoothing estimates for scalar
 equations]\label{prop:sm}
  Let $n\geq2$. Assume the operator
  $$-\Delta+W(x,D)=-\Delta+a(x)\cdot \nabla+b_{1}(x)+b_{2}(x)$$
  is selfadjoint with
  \begin{equation}\label{eq.ipab}
    |a(x)|\le \frac{\delta}{\tau_{\epsilon}w_{\sigma}^{1/2}},\qquad
    |b_{1}(x)|\le \frac{\delta}{\tau_{\epsilon}^{2}},\qquad
    0\le b_{2}(x)\le \frac{C}{\tau_{\epsilon}^{2}}
  \end{equation}
  for some $\delta,\epsilon>0$ sufficiently small and
  some $\sigma>1$, $C>0$. Moreover assume that
  0 is not a resonance for $-\Delta+b_{2}$.

  Then the following smoothing estimates hold: for the
  Schr\"odinger equation
   \begin{equation*}
     \|\tau_\epsilon^{-1}e^{it(-\Delta+W)}f\|_{L^2L^2}+
     \|\tau_\epsilon^{-1}|D|^{1/2}e^{it(-\Delta+W)}f\|_{L^2L^2}
   \lesssim\|f\|_{L^2}
  \end{equation*}
  while for the wave and Klein-Gordon equations
   \begin{equation*}
     \|\tau_\epsilon^{-1}e^{it\sqrt{-\Delta+W}}f\|_{L^2L^2}+
     \|\tau_\epsilon^{-1}e^{it\sqrt{1-\Delta+W}}f\|_{L^2L^2}
   \lesssim\|f\|_{L^2}.
  \end{equation*}
\end{proposition}

The assumption that 0 is not a resonance for $-\Delta+b_{2}(x)$
here means: if $(-\Delta+b_{2})f=0$ and $\xx^{-1}f\in L^{2}$ then
$f \equiv0$.

We can  then prove
Strichartz estimates for the perturbed scalar equations
as a consequence of the above smoothing properties.
Notice that we must require some additional
regularity on the magnetic coefficient $a(x)$.
Moreover, the use of the Christ-Kiselev lemma
(see Section \ref{sec.stri} for details) prevents us from reaching
the endpoint.

\begin{theorem}[Strichartz for Schr\"odinger]\label{thm.schrogrande}
  Let $n\ge2$, $-\Delta+W$ be as in Proposition \ref{prop:sm}
  and assume in addition that
  \begin{equation}\label{eq.ipab2}
    \xx^{1+3\epsilon}\chi(x)a_j(x)\in C^{\frac12+2\epsilon}  \qquad
    \text{for some function}\quad \chi\gtrsim w_{\sigma}^{1/2}.
  \end{equation}
  Then, for any non-endpoint Schr\"odinger admissible
  couple $(p,q)$, the following Strichartz estimate holds:
  \begin{equation}\label{eq.stri.schro}
    \|e^{it(-\Delta+W)}f\|_{L^pL^q}\lesssim\|f\|_{L^2}.
  \end{equation}
\end{theorem}

\begin{theorem}[Strichartz for wave]\label{thm.onde}
  Let $n\ge3$, $-\Delta+W$ be as in Proposition \ref{prop:sm}
  and assume in addition that
  \begin{equation}\label{eq.ipabonde}
    |a(x)|\le \frac{C}{\tau_{\epsilon}^{2}},\qquad
    |b_{1}+b_{2}-\nabla \cdot a|\le
     \frac{C}{|x|\tau_{\epsilon}}.
  \end{equation}
  Then, for any non-endpoint wave admissible couple $(p,q)$ the
  following Strichartz estimate holds:
  \begin{equation}\label{eq.strifinonde}
    \|e^{it\sqrt{-\Delta+W}}f\|_{L^p\dot H^{\frac1q-\frac1p-\frac12}_q}
    \lesssim\|f\|_{L^{2}}.
  \end{equation}
\end{theorem}

\begin{theorem}[Strichartz for Klein-Gordon]\label{thm.kleingordon}
  Let $n\ge2$, $-\Delta+W$ be as in Proposition \ref{prop:sm}
  and assume in addition that
  \begin{equation}
    |a(x)|\le \frac{C}{\tau_{\epsilon}^{2}},\qquad
    |b_{1}+b_{2}-\nabla \cdot a|\le
     \frac{C}{\xx\tau_{\epsilon}}.
  \end{equation}
  Then, for any non-endpoint Schr\"odinger admissible
  couple $(p,q)$, the following Strichartz estimate holds:
  \begin{equation}
    \|e^{it\sqrt{-\Delta+1+W}}\|_{L^p H^{\frac1q-\frac1p-\frac12}_q}\leq
    C\|f\|_{L^{2}}.
  \end{equation}
\end{theorem}

Our final results concern the Dirac system:

\begin{theorem}[Massless Dirac]\label{thm.dirac}
  Let $n=3$, and let $V(x)=V(x)^{*}$ be a
  4$\times$4 complex valued matrix such that
  \begin{equation}\label{eq.ipVfin}
    |V(x)|\leq\frac{\delta}{w_{\sigma}(x)}
  \end{equation}
  for some $\delta$ sufficiently small and some $\sigma>1$.
  Then the following smoothing estimate holds:
  \begin{equation}\label{eq:smdi}
    \|w_{\sigma}^{-1/2}e^{it(\D+V)}f\|_{L^{2}L^{2}}\lesssim
      \|f\|_{L^{2}}
  \end{equation}
  and, for any non-endpoint wave admissible couple $(p,q)$, the
  following Strichartz estimate holds:
  \begin{equation}\label{eq.stridir1}
    \|e^{it(\D+V)}f\|_{L^p\dot H^{\frac1q-\frac1p-\frac12}}\lesssim
    \|f\|_{L^{2}}.
  \end{equation}
\end{theorem}

\begin{theorem}[Massive Dirac]\label{thm.diracmass}
  Let $n=3$, and let $V(x)=V(x)^{*}$ be a
  4$\times$4 complex valued matrix such that
  \begin{equation}\label{eq.ipVfin2}
    |V(x)|\leq\frac{\delta}{\tau_{\epsilon}(x)}
  \end{equation}
  for some $\delta,\epsilon>0$ sufficiently small.
  Then the following smoothing estimate holds:
  \begin{equation}\label{eq:smdim}
    \|\tau_{\epsilon}^{-1} e^{it(\D+\beta+ V)}f\|_{L^{2}L^{2}}\lesssim
      \|f\|_{L^{2}}
  \end{equation}
  and, for any non-endpoint
  Schr\"odinger admissible couple $(p,q)$, the
  following Strich\-artz estimate holds:
  \begin{equation}\label{eq.stridir2}
    \|e^{it(\D+\beta+ V)}f\|
    _{L^p H^{\frac1q-\frac1p-\frac12}}\lesssim
    \|f\|_{L^{2}}.
  \end{equation}
\end{theorem}

The paper is organized as follows: in Section 
\ref{sec.resolv} we prove resolvent estimates for the
perturbed operator, which are equivalent to smoothing
estimates for the corresponding flow via 
Kato theory, while Section \ref{sec.stri} is devoted
to the proof of the main theorems. A short Appendix
collects the estimates for the free
Klein-Gordon and Dirac equations; these can be obtained
by a standard application of the Ginibre-Velo and
Keel-Tao methods, and we decided to include a sketch
of the proof for the sake of completeness.

\section{Resolvent Estimates}\label{sec.resolv}

In this section we shall prove the basic resolvent estimates
for the perturbed operators, which are the crucial step in
the proof. As an immediate consequence we shall obtain
smoothing estimates for the corresponding evolution
operators, by a standard application of the well-known result of
Kato (see \cite{reed-simonII}):

\begin{theorem}[Kato smoothing Theorem, \cite{Kato66}]
  \label{thm.kato}
    Let $X,Y$ be Hilbert spaces, let $H:X\to X$
    be a self-adjoint operator whose resolvent we denote by
    $R(\lambda)=(H-\lambda)^{-1}$,
    and let $A:X\to Y$ be a closed,
    densely defined operator, which may be unbounded. Assume that
    \begin{equation}\label{eq.kato}
      \|AR(\lambda)A^*g\|_Y\le M \|g\|_{Y}\qquad
         \forall g\in D(A^*),\ \lambda\not\in \mathbb{R}.
    \end{equation}
    Then the operator $A$ is $H$-smooth, i.e.,
    $e^{itH}f\in D(A)$ for all $f\in X$ and a.e. $t$, and
    \begin{equation}\label{eq.tesikato}
        \int_{-\infty}^\infty
        \|Ae^{-itH}f\|^{2}_Y\,dt\leq\frac{2}{\pi}M^{2}\|f\|^{2}_{X}\qquad
           \forall f\in X.
    \end{equation}
\end{theorem}

\subsection{The magnetic Schr\"odinger operator}\label{ssec.reslap}

The following lemma gives sufficient conditions for the
magnetic Schr\"odinger operator
$H=-(\nabla+iA(x))^{2}+B(x)$ to be selfadjoint.
We sketch a proof since the assumptions on the
coefficients are not completely standard:

\begin{lemma}\label{lem.selfadj}
   Let $A_{j}(x)$, $A=(A_{1},\dots,A_{n})$
   and $B(x)$ be real valued functions satisfying
   \begin{equation}\label{eq.assco}
       \|B_{+}\|_{L^{n/2,\infty}}<C,\qquad
       \|B_{-}\|_{L^{n/2,\infty}}<\delta,\qquad
       \|A\|_{L^{n,\infty}}<\delta
   \end{equation}
   for some $C,\delta>0$.
   Then, if $\delta$ is sufficiently small, the operator
   \begin{equation}\label{eq.operDW2}
     H=-(\nabla+iA(x))^{2}+B(x)
   \end{equation}
   can be uniquely defined as a selfadjoint nonnegative
   operator in $L^{2}$, with form domain $H^{1}(\mathbb{R}^{n})$.
   Moreover we have
   \begin{equation}\label{eq.posit}
    \|H^{1/2}g\|_{L^{2}}\simeq\|g\|_{\dot H^{1}}.
   \end{equation}
\end{lemma}

\begin{proof}
  The quadratic form
  \begin{equation*}
    q(\phi,\psi)=
    \left((\nabla+iA(x))\phi,(\nabla+iA(x))\psi\right)_{L^{2}}
        +(B(x)\phi,\psi)_{L^{2}}
  \end{equation*}
  is well defined on $H^1 \times H^{1}$ under assumptions
  \eqref{eq.assco}. Indeed, using
  the embedding $H^1\subset L^{2n/(n-2),2}$,
  H\"older's inequality in Lorentz spaces \cite{oneil}
  and assumptions \eqref{eq.assco}, we have
  \begin{align*}
    |q(\varphi,\varphi)| \le&
    \|\nabla\varphi\|_{L^2}^2+
    2\|A\|_{L^{n,\infty}}\|\nabla\varphi\cdot
    \overline\varphi\|_{L^{\frac{n}{n-1},1}}+
    \||A|^{2}+|B|\|_{L^{\frac n2,\infty}}\|\varphi^{2}
       \|_{L^{\frac{n}{n-2},1}}
       \\
    \lesssim&\|\nabla\varphi\|_{L^2}^2.
  \end{align*}
The form $q$ is symmetric since $A$ and $B$ are real
valued. By standard results (see e.g.~\cite{reed-simonII},
Theorem VIII.15), $q$ is the form associated to a
unique defined self-adjoint operator provided
the form is {\it closed},
i.e. its domain $H^1(\R^n)$ is complete under the norm
\begin{equation}\label{eq.normeq}
  \triple{\varphi}^2=q(\varphi,\varphi)
  +C\|\varphi\|_{L^2}^2,
\end{equation}
for some $C>0$, and it is {\it semibounded}, i.e.
\begin{equation}\label{eq.semib}
  q(\varphi,\varphi)\geq-C\|\varphi\|_{L^2}^2,
\end{equation}
for some $C>0$. To prove this we estimate the form
from below as follows
\begin{align*}
  q(\varphi,\varphi) &=
  \|\nabla\varphi\|_{L^2}^2
    +2\Im(A\cdot\nabla\varphi,\varphi)_{L^2}
    +((|A|^{2}+B_+)\varphi,\varphi)_{L^2}
    -(B_-\varphi,\varphi)_{L^2}
    \\ &\geq
    \|\nabla\varphi\|_{L^2}^2
    +2\Im(A\cdot\nabla\varphi,\varphi)_{L^2}
    -(B_-\varphi,\varphi)_{L^2}.
\end{align*}
Proceeding as for the upper bound we obtain
\begin{equation}\label{eq.stimabasso}
  q(\varphi,\varphi)\geq
  \|\nabla\varphi\|_{L^2}^2
  -C\delta\|\nabla\varphi\|_{L^2}^2
  \gtrsim
  \|\nabla\varphi\|_{L^{2}}^{2}
\end{equation}
for $\delta$ small enough.
This proves the semiboundedness of the form and
\eqref{eq.posit}, which implies that the norm
\eqref{eq.normeq} is equivalent to the norm
of $H^{1}$ and hence the form is closed.
\end{proof}

%

We now investigate in some detail the properties
of the resolvent operators
\begin{equation}\label{eq.RW}
  R(z)=(-\Delta+W-z)^{-1}
\end{equation}
\begin{equation*}
  R_0(z)=(-\Delta-z)^{-1},\qquad
  R_{b_{2}}(z)=(-\Delta+b_{2}(x)-z)^{-1}.
\end{equation*}
The following weight functions will appear in our
resolvent estimates ($\epsilon>0,\sigma>1$):
\begin{equation}\label{eq.weights}
  \xx=(1+|x|^2)^{\frac12},\qquad
  w_\sigma(x)=|x|(1+|\log|x||)^\sigma,
\end{equation}
and
\begin{equation}\label{eq.weights2}
  \tau_\epsilon(x)=
  \begin{cases}
      |x|^{\frac12-\epsilon}+|x|&\text{if $n\ge3$,}\\
      |x|^{\frac12-\epsilon}+|x|^{1+\epsilon}&\text{if $n=2$.}
  \end{cases}
\end{equation}
Notice that
\begin{equation*}
  |x|\leq\tau_\epsilon(x),\qquad
  w_\sigma^{\frac12}(x)\leq C\tau_\epsilon(x),\qquad
\end{equation*}
and
\begin{equation*}
  \tau_\epsilon(x)\leq C\xx,\text{ for $n\ge3$},\qquad
  \tau_\epsilon(x)\leq C\xx^{1+\epsilon},\text{ for $n=2$}
\end{equation*}
for some constant $C=C(\epsilon,\sigma)$.

In order to estimate the resolvent $R$ we shall use
the formal identity
\begin{equation}\label{eq.idreslap}
  R=R_0(I+b_{2}R_0)^{-1}
  (I+(a\cdot\nabla +b_{1})R_{b_{2}})^{-1}.
\end{equation}
Our first goal
will be to prove that the operators
$(I+bR_0)^{-1}$
and $(I+(a\cdot\nabla +b_{1})R_{b_{2}})^{-1}$
are well defined and uniformly bounded
in suitable weighted $L^{2}$ spaces.
In the following lemma, the assumption that
0 is not a resonance of $-\Delta+b(x)$ means that the only
distribution solution $f$ of the equation
$-\Delta f+bf=$ belonging to $L^{2}(\xx^{-2}dx)$ is
$f \equiv0$.

\begin{lemma}\label{lem.invers1}
  Let $b(x)$ be real valued and such that, for some
  $\epsilon,\delta>0$ small enough (recall \eqref{eq.weights2}),
  \begin{equation}\label{eq.ipinvers1}
    \|\tau_\epsilon^2b_{+}\|_{L^\infty}<\infty,\qquad
    \|\tau_\epsilon^2b_{-}\|_{L^\infty}<\delta.
  \end{equation}
  Assume that 0 is not a resonance for $-\Delta+b(x)$.
  Then $I+bR_0(z)$ is invertible with a uniformly
  bounded inverse on $L^2(\tau_\epsilon^2dx)$:
  \begin{equation}\label{eq.estinvers1}
    \|\tau_\epsilon(I+bR_0(z))^{-1}f\|_{L^2}\leq
    C\|\tau_\epsilon f\|_{L^2}.
  \end{equation}
\end{lemma}

\begin{proof}
We recall the following estimates for
the free resolvent $R_0$: fix any $\sigma>1$, then
for all $z\in\C$
\begin{equation}\label{eq.brv}
  \|w_\sigma^{-\frac12}R_0(z)f\|_{L^2}\leq
  \frac{C}{\sqrt{|z|}}\|w_\sigma^{\frac12}f\|_{L^2},
\end{equation}
\begin{equation}\label{eq.brv2}
  \|w_\sigma^{-\frac12}\nabla R_0(z)f\|_{L^2}\leq
  C\|w_\sigma^{\frac12}f\|_{L^2},
\end{equation}
\begin{equation}\label{eq.KY1}
  \|\,|x|^{-1}R_0f\|_{L^2}\leq C\|\,|x|f\|_{L^2},\qquad
  n\ge3
\end{equation}
\begin{equation}\label{eq.KY2}
  \|\,|x|^{-1+\epsilon}|D|^{\epsilon}R_0f\|_{L^2}\leq
  C\|\,|x|^{1-\epsilon}|D|^{-\epsilon}f\|_{L^2},\qquad
  n=2\quad(0<\epsilon<1/2)
\end{equation}
(see \cite{BRV}, \cite{DAnconaFanelli06-c}
for \eqref{eq.brv}, \eqref{eq.brv2}, and
\cite{KY} for \eqref{eq.KY1}-\eqref{eq.KY2}).
As usual, for $\lambda\in\R^{+}$ the resolvent $R_{0}(z)$
must be replaced with
the limit operators $R_{0}(\lambda\pm i0)$.
By the elementary inequalities
$|x|\leq\tau_\epsilon(x)$, $w_\sigma^{\frac12}(x)
\leq C\tau_\epsilon(x)$, we can condense the estimates
\eqref{eq.brv} and \eqref{eq.KY1} in the following (weaker) one
for $n\ge3$:
\begin{equation}\label{eq.tau}
  \|\tau_\epsilon^{-1}R_0(z)f\|_{L^2}\leq
  \frac C{\sqrt{\langle z\rangle}}
  \|\tau_\epsilon f\|_{L^2},
  \qquad\text{for all}\ z\in\C.
\end{equation}
In dimension $n=2$ we deduce by duality from \eqref{eq.KY2}
the following
\begin{equation*}
  \|\,|D|^{\epsilon}|x|^{-1+\epsilon}R_0f\|_{L^2}\leq
  C\|\,|D|^{-\epsilon}|x|^{1-\epsilon}f\|_{L^2},\qquad
\end{equation*}
which implies, via Sobolev embedding and H\"older inequality,
\begin{equation*}
  \|\xx ^{-\sigma}|x|^{-1+\epsilon}R_0f\|_{L^2}\leq
  C\|\xx ^{\sigma}|x|^{1-\epsilon}f\|_{L^2},\qquad
  \sigma>\epsilon
\end{equation*}
and hence \eqref{eq.tau} follows also for $n=2$
(recall \eqref{eq.weights2})

Now, using assumption \eqref{eq.ipinvers1}, we have
\begin{equation}\label{eq.neum}
  \|\tau_\epsilon bR_0(z)f\|_{L^2}\leq
  \|\tau_\epsilon^2b\|_{L^\infty}
  \|\tau_\epsilon^{-1}R_0(z)f\|_{L^2}\leq
  \frac C{\sqrt{\langle z\rangle}}
  \|\tau_\epsilon^2b\|_{L^\infty}
  \|\tau_\epsilon f\|_{L^2},
\end{equation}
with $C$ as in \eqref{eq.tau}; hence, if $z$ is
sufficiently large, namely so large that
\begin{equation*}
  \langle z\rangle>
  C^2\|\tau_\epsilon^2b\|_{L^\infty}^2,
\end{equation*}
we can invert the operator $I+bR_0$ by a Neumann series
in the weighted space $L^{2}(\tau_{\epsilon}^{2}dx)$,
with a uniform bound on the norm of the inverse.

In the low frequency case
\begin{equation}\label{eq.region}
  \langle z\rangle\leq
  C^2\|\tau_\epsilon^2b\|_{L^\infty}^2,
\end{equation}
the family of operators $(I+bR_0)(z)$ is uniformly bounded in
$L^2(\tau_\epsilon^2dx)$ by
\eqref{eq.neum}. We also notice that $bR_0$ is a compact operator
on $L^2(\tau_\epsilon^2dx)$; indeed, $R_{0}$ is a compact operator
from $L^{2}(\tau_{\epsilon}^{2}dx)$ to $L^{2}(\tau_{\epsilon}^{-2}dx)$
(see \eqref{eq.brv}--\eqref{eq.brv2}), while multiplication
by $b$ is  bounded
from $L^{2}(\tau_{\epsilon}^{-2}dx)$ to $L^{2}(\tau_{\epsilon}^{2}dx)$.
Thus by standard analytic Fredholm theory
we can invert $I+bR_0(z)$ uniformly in $z$,
provided $I+bR_0(z)$ is injective
on $L^2(\tau_\epsilon^2dx)$ for each fixed $z$.
This is obvious for $z$ outside $\overline{\mathbb{R}^{+}}$,
since by our assumptions the operator $-\Delta+b$
is nonnegative and selfadjoint, and
is true by assumption for $z=0$, hence we need only check
the case $z=\lambda>0$.

Thus let $\lambda\ge0$ and $f\in L^2(\tau_\epsilon^2dx)$
such that $f+b(x)R_{0}(\lambda+i0)f=0$
(the $-i0$ case is identical). We notice that
estimate \eqref{eq.brv2} implies that
$R_{0}(z)f\in H^{1}_{\text{loc}}$ and hence in particular
$R_{0}(z)f$ is in $L^{2n/(n-1)}$ locally. Since
$|b|\lesssim \tau_{\epsilon}^{-2}$ which is locally in
$L^{n}$, we conclude that $f=-bR_{0}(\lambda)f$ is locally
in $L^{2}$. Recalling that $f\in L^2(\tau_\epsilon^2dx)$
this implies $f\in L^2(\xx^{2}dx)$. Thus we are in the
framework of the standard Agmon theory and we deduce that
$\lambda$ is an eigenvalue of $-\Delta+b(x)$; but this
is excluded under our assummtions on $b$, for instance by the
results of\cite{ionescujerison} (Theorem 2.1).

In conclusion, we can invert
$(I+bR_0)(z)$ in $L^2(\tau_\epsilon^2dx)$
with an uniform bound for the inverse $(I+bR_0)^{-1}$,
and this completes the proof.
\end{proof}

The preceding lemma allows us to construct the
resolvent operator
\begin{equation}\label{eq:firstres}
  R_{b}(z)=R_{0}(z)(I+bR_{0}(z))^{-1},
\end{equation}
which, in view of \eqref{eq.estinvers1} and
\eqref{eq.tau}, is a bounded operator from $L^{2}(\tau_{\epsilon}^{2}dx)$
to $L^{2}(\tau_{\epsilon}^{-2}dx)$ for all $z \in \mathbb{C}$.

We have next:

\begin{lemma}\label{lem.invers2}
  Consider the operator $-\Delta+a(x)\cdot \nabla+b_{1}(x)+b_{2}(x)$
  under the following assumptions: the operator is selfadjoint,
  $b_{2}$ is real valued and nonnegative, and for some
  $\delta,\epsilon>0$ small enough, $\sigma>1$,
  \begin{equation}\label{eq.ipinvers2}
    \|\tau_\epsilon w_\sigma^{\frac12}
    a\|_{L^\infty}+
    \|\tau_\epsilon^2b_{1}\|_{L^\infty}< \delta,\qquad
    \|\tau_\epsilon^2b_{2}\|_{L^\infty}<\infty.
  \end{equation}
  Moreover assume that 0 is not a resonance for
  $-\Delta+b_{2}(x)$. Then
  $I+(a\cdot\nabla+b_{1})R_{b_{2}}$ is invertible with a bounded inverse
  on $L^2(\tau_\epsilon^2dx)$:
  \begin{equation}\label{eq.estinvers2}
    \|\tau_\epsilon(I+(a\cdot\nabla+ b_{1})R_{b_{2}})^{-1}
    f\|_{L^2}\leq C\|\tau_\epsilon f\|_{L^2}.
  \end{equation}
\end{lemma}

\begin{proof}
Using assumptions \eqref{eq.ipinvers2}, H\"older inequality
and estimate \eqref{eq.brv2}, we can write
\begin{eqnarray*}
  \|\tau_\epsilon a\cdot\nabla R_{b_{2}}f\|_{L^2}
  & \leq & \|\tau_\epsilon a\cdot\nabla R_0
  (I+b_{2}R_0)^{-1}f\|_{L^2} \\ \  & \leq &
  \|\tau_\epsilon w_\sigma^{\frac12}
  a\|_{L^\infty}\|w_\sigma^{-\frac12}
  \nabla R_0(I+b_{2}R_0)^{-1}f\|_{L^2} \\ \  & \lesssim &
  \delta\cdot\|w_\sigma^{\frac12}(I+bR_0)^{-1}f\|_{L^2}
  \\ \  & \lesssim &
  \delta\cdot\|\tau_\epsilon(I+b_{2}R_0)^{-1}f\|_{L^2}
\end{eqnarray*}
and Lemma \ref{lem.invers1} gives finally
\begin{equation*}
  \|\tau_\epsilon a\cdot\nabla R_{b_{2}}f\|_{L^2}
  \lesssim\delta\cdot \|\tau_\epsilon f\|_{L^2}.
\end{equation*}
On the other hand, by \eqref{eq.ipinvers2}
and estimate \eqref{eq.tau}
\begin{align*} 
  \|\tau_\epsilon b_{1} R_{b_{2}}f\|_{L^2}  \leq &\
     \|\tau_\epsilon^{2}b_{1}\|_{L^\infty}
     \|\tau_\epsilon^{-1} R_0(I+b_{2}R_0)^{-1}f\|_{L^2} \\
  \lesssim &\ \delta \cdot
     \|\tau_\epsilon(I+b_{2}R_0)^{-1}f\|_{L^2}
\end{align*}
and again by Lemma \ref{lem.invers1} we have
\begin{equation*}
  \|\tau_\epsilon b_{1} R_{b_{2}}f\|_{L^2}
  \lesssim\delta\cdot \|\tau_\epsilon f\|_{L^2}.
\end{equation*}
Thus, if $\delta$ is sufficiently small, we can invert
$I+(a\cdot \nabla+b_{1}) R_{b_{2}}$
via a Neumann series, and we obtain \eqref{eq.estinvers2}.
\end{proof}

We collect and complete the above estimates in the
following

\begin{proposition}\label{prop.estreslap}
  Consider the operator
  $-\Delta+W(x,D)\equiv-\Delta+a(x)\cdot \nabla+b_{1}(x)+b_{2}(x)$
  under the assumptions: the operator is selfadjoint,
  $b_{2}$ is real valued and nonnegative, and for some
  $\delta,\epsilon>0$ small enough, $\sigma>1$,
  \begin{equation}\label{eq.ab}
    \|\tau_\epsilon w_\sigma^{\frac12}
    a\|_{L^\infty}+
    \|\tau_\epsilon^2b_{1}\|_{L^\infty}< \delta,\qquad
    \|\tau_\epsilon^2b_{2}\|_{L^\infty}<\infty.
  \end{equation}
  Moreover assume that 0 is not a resonance for
  $-\Delta+b_{2}(x)$.
  Then the resolvent operator  $R(z)=(-\Delta+W-z)^{-1}$
  satisfies the following estimates for all $z\in \mathbb{C}$:
  \begin{equation}\label{eq.estRW}
    \|\tau_\epsilon^{-1}R(z)f\|_{L^2}\leq
    \frac C{\sqrt{\langle z\rangle}}\|
    \tau_\epsilon f\|_{L^2},
  \end{equation}
  \begin{equation}\label{eq.estRW2}
    \|\tau_\epsilon^{-1}\nabla
    R(z)f\|_{L^2}\leq
    C\|\tau_\epsilon f\|_{L^2}.
  \end{equation}
  and
  \begin{equation}\label{eq.estlapfinal}
    \|\xx^{-1}R(z)f\|_{H^1}\leq C\|\xx f\|_{L^2},\qquad
    n\ge3;
  \end{equation}
  replace the weights $\xx^{-1},\xx$ by
  $\xx^{-1-\epsilon},\xx^{1+\epsilon}$ respectively
  in dimension 2.
  As a consequence, the Schr\"odinger flow
  $e^{it(-\Delta+W)}f$ has the
  smoothing property
  \begin{equation}\label{eq.katoschro1}
    \|\tau_\epsilon^{-1}e^{it(-\Delta+W)}f\|_{L^2L^2}+
    \|\tau_\epsilon^{-1}|D|^{1/2}e^{it(-\Delta+W)}f\|_{L^2L^2}
  \leq C\|f\|_{L^2}.
 \end{equation}
\end{proposition}

\begin{remark}\label{rem:smoo}
  For the following applications it will be convenient to rewrite
  the (second) smoothing estimate above in the equivalent form
   \begin{equation}\label{eq.katoschro1b}
     \|\tau_\epsilon^{-1}{\nabla}{|D|^{-1/2}}
        e^{it(-\Delta+W)}f\|_{L^2L^2}
   \leq C\|f\|_{L^2}.
  \end{equation}
  This follows immediately from the fact that
  $\partial_{j}|D|^{-1/2}=iR_{j}|D|^{1/2}$, where
  $R_{j}=i^{-1}\partial_{j}|D|^{-1}$ is the $j$-th Riesz
  operator, and on the other hand $\tau_{\epsilon}^{-1}$
  is an $A_{2}$ weight, as proved in Lemma \ref{lem.A2} below.
\end{remark}

\begin{proof}
Estimates \eqref{eq.estRW} and \eqref{eq.estRW2}
are immediate conswequences of \eqref{eq.idreslap},
\eqref{eq.brv2} and of Lemmas
\ref{lem.invers1}, \ref{lem.invers2}. Moreover,
\eqref{eq.estRW} implies in particular
\begin{equation*}
    \|\tau_\epsilon^{-1}R(z)f\|_{L^2}\leq C
    \|\tau_\epsilon f\|_{L^2},
\end{equation*}
and the Kato smoothing theorem with the choices
$A=\tau_{\epsilon}^{-1}$, $X=Y=L^{2}$ gives the
first estimate in \eqref{eq.katoschro1}.

To prove \eqref{eq.estlapfinal}, write
\begin{equation*}
  \|\xx^{-1}Rf\|_{H^1}\lesssim
  \|\xx^{-1}Rf\|_{L^2}
  +\|\xx^{-2}Rf\|_{L^2}+\|\xx^{-1}
  \nabla Rf\|_{L^2}.
\end{equation*}
The first two terms can be estimated by \eqref{eq.estRW}
\begin{equation}\label{eq.I}
  \|\xx^{-2}Rf\|_{L^2}+
  \|\xx^{-1}Rf\|_{L^2}\leq
  C\|\tau_{\epsilon}^{-1} Rf\|_{L^2}\leq
  C\|\tau_\epsilon f\|_{L^2}\leq
  C\|\xx f\|_{L^2},
\end{equation}
while the third term is bounded using \eqref{eq.estRW2}:
\begin{equation}\label{eq.III}
  \|\xx^{-1}\nabla Rf\|_{L^2}\leq
  \|\tau_\epsilon^{-1} Rf\|_{L^2}\leq
  C\|\tau_\epsilon f\|_{L^2}\leq
  C\|\xx f\|_{L^2}
\end{equation}
and this proves \eqref{eq.estlapfinal}.

Now write \eqref{eq.estlapfinal} in the equivalent forms
\begin{equation*}
   \|\langle D\rangle\xx^{-1}R(z)\xx^{-1}f\|_{L^{2}}
    \leq C\| f\|_{L^2}
\end{equation*}
and, by duality,
\begin{equation*}
   \|\xx^{-1}R(z)\xx^{-1}\langle D\rangle f\|_{L^{2}}
   \leq C\|f\|_{L^2}.
\end{equation*}
The last two estimates state that the operator
$\xx^{-1} R(z)\xx^{-1}$ is bounded, uniformly in $z\in\C$,
from $L^{2}$ to $H^{1}$ and from $H^{-1}$ to $L^{2}$.
By complex interpolation this implies that it is also
bounded from $H^{-1/2}$ to $H^{1/2}$, i.e.,
\begin{equation*}
   \|\langle D\rangle^{1/2}
    \xx^{-1}R(z)\xx^{-1}\langle D\rangle^{1/2}
     f\|_{L^{2}}
   \leq C\|f\|_{L^2}
\end{equation*}
Then by Kato smoothing we obtain
also the second estimate in \eqref{eq.katoschro1}.

The proof for the case $n=2$ is completely analogous.
\end{proof}

\subsection{The wave and Klein-Gordon generators.}

We consider now the operator $\sqrt{-\Delta+W}$,
where as usual
$$W=W(x,D)=a\cdot\nabla+b,\qquad b=b_{1}+b_{2}$$
which generates the flow $e^{it \sqrt{-\Delta+W}}$
of the perturbed wave equation. The free operator
$|D|:=\sqrt{-\Delta}$ is self-adjoint and nonnegative
on $L^2$, and can be handled as follows. If we denote its
resolvent by $R_{|D|}(z)=(|D|-z)^{-1}$, we have
\begin{equation}\label{eq.idresD}
  R_{|D|}(z)=(|D|+z)R_0(z^2).
\end{equation}
This simple identity allows us to estimate $R_{|D|}$
using some standard
techniques from harmonic analysis. We need a lemma:

\begin{lemma}\label{lem.A2}
    Let $n\geq2$.
    For any $\sigma>1$, the weight
    $w_{\sigma}=|x|(1+|\log|x||)^{\sigma}$ is an
    $A_{2}$ weight, i.e., there exist a constant $A$ such
    that, for any ball $B=B(x_{0},R)$,
    \begin{equation}\label{eq.A2}
     A(x_{0},R)\equiv
      \left[\frac{1}{|B|}\int_{B}w_{\sigma}dx\right]\cdot
      \left[\frac{1}{|B|}\int_{B}w_{\sigma}^{-1}dx\right]\leq
         A<\infty.
    \end{equation}
    Obviously, we have also $w_{\sigma}^{-1}\in A_{2}$.
    The same property holds for the weights $\tau_{\epsilon}$,
    $\tau_{\epsilon}^{-1}$ defined in \eqref{eq.weights2}.
\end{lemma}

\begin{proof}
The bound for the function $A(x_{0},R)$ is trivial
if $R\le|x_{0}|/2$, indeed it is sufficient
to write
\begin{equation*}
    A(x_{0},R)\leq
    C\max_{B}w_{\sigma}\cdot
      \max_{B}w_{\sigma}^{-1}
      \le C'
\end{equation*}
since the ball $B$ is at a distance greater than $|x_{0}|/2$ from the
origin.

If, on the other hand, $R\ge|x_{0}|/2$, it is easy to check
that $A(x_{0},R)$ is bounded by a constant (depending only
on the space dimension $n$) times $A(0,3R)$. Thus we are
reduced to the case of balls $B(0,R)$ centered in 0.

For \emph{small} $R\le 10$
the function $A(0,R)$ is bounded. Indeed, H\^opital's
theorem gives
\begin{equation*}
    \lim_{\epsilon\downarrow0}\int_{0}^{\epsilon}
       \frac{r^{n-2}dr}{(1+|\log r|)^{\sigma}}
       \cdot
       \frac{(1+|\log\epsilon|)^{\sigma}}
           {\epsilon^{n-1}}=\frac 1{n-1}
\end{equation*}
which implies for small $R$
\begin{equation}\label{eq.smallR}
    \int_{0}^{R}
       \frac{r^{n-2}dr}{(1+|\log r|)^{\sigma}}\sim
       \frac{R^{n-1}}{(1+|\log R|)^{\sigma}}
\end{equation}
and similarly
\begin{equation*}
    \int_{0}^{R}
       r^{n}(1+|\log r|)^{\sigma}dr\sim
       R^{n+1}(1+|\log R|)^{\sigma}
\end{equation*}
whence we get $A(0,R)\le C$.

For \emph{large} $R>10$ we rescale and obtain
\begin{equation*}
    A(0,R)=\int_{0}^{1}
      \frac{\tau^{n-2}d\tau}{(1+|\log R+\log\tau|)^{\sigma}}
      \cdot\int_{0}^{1}
      \tau^{n}(1+|\log R+\log\tau|)^{\sigma}d\tau
\end{equation*}
The second integral is clearly bounded by $C(\log R)^{\sigma}$.
The first integral can be split into
\begin{equation*}
    \int_{0}^{1/\sqrt R}
      \frac{\tau^{n-2}d\tau}{(1+|\log R+\log\tau|)^{\sigma}}\le
    \int_{0}^{1/\sqrt R}
      \frac{\tau^{n-2}d\tau}{(1+|\log\tau|)^{\sigma}}\sim
      \frac{R^{-\frac{n-1}2}}{(1+\frac12\log R)^{\sigma}}\le
      R^{-\frac12}
\end{equation*}
where we used again \eqref{eq.smallR}, and
\begin{equation*}
    \int_{1/\sqrt R}^{1}
      \frac{\tau^{n-2}d\tau}{(1+|\log R+\log\tau|)^{\sigma}}\le
    \int_{1/\sqrt R}^{1}
      \frac{\tau^{n-2}d\tau}{(1+\frac12\log R)^{\sigma}}\le
    C(\log R)^{-\sigma}.
\end{equation*}
Putting everything together, we obtain the required bound
also for large $R$, and this concludes the proof of the Lemma.

The proof for $\tau_{\epsilon}$ is much simpler. We reduce as
above to the case of spheres $B(0,R)$ centered in the origin.
For $R\le1$ we can use the equivalence
$\tau_{\epsilon}\simeq |x|^{1/2-\epsilon}$ and the bound follows
from the well-known fact that $|x|^{1/2-\epsilon}$ is an
$A_{2}$ weight. For $R>1$ we use the estimate
\begin{equation*}
  A(0,R)\lesssim
   \frac{1}{|B|}\int_{B}(1+|x|)dx \cdot
   \frac{1}{|B|}\int_{B}\frac{dx}{|x|}
\end{equation*}
(replace $|x|$ with $|x|^{1+\epsilon}$ for $n=2$)
whence the bound follows easily.
\end{proof}

Knowing that $w_{\sigma}^{-1}\in A_{2}$, we see that
the Riesz operators
\begin{equation*}
    R_{j}=i^{-1}\frac{\partial_{j}}{|D|}
\end{equation*}
are bounded on the space $L^{2}(w_{\sigma}^{-1}dx)$
by standard results
(see e.g. the Corollary to Theorem 2, \S V.4.2
of \cite{St}). Writing $|D|=i\sum R_{j}\partial_{j}$,
we have
\begin{equation*}
    \|w_{\sigma}^{-1/2}|D|g\|_{L^{2}}\le
    \sum_{j}\|w_{\sigma}^{-1/2}R_{j}\partial_{j}g\|_{L^{2}}\le
    C\|w_{\sigma}^{-1/2}\nabla g\|_{L^{2}}.
\end{equation*}
Thus estimate \eqref{eq.idresD} implies
\begin{equation}\label{eq.res1}
    \|w_\sigma^{-\frac12}R_{|D|}(z)f\|_{L^2}\le
    C\|w_\sigma^{-\frac12}\nabla R_{0}(z)f\|_{L^{2}}+
    C|z|\cdot\|w_\sigma^{-\frac12} R_{0}(z)f\|_{L^{2}}.
\end{equation}
Then, inequalities \eqref{eq.brv} and \eqref{eq.brv2} yield
immediately the following estimate for the
free resolvent: for any fixed $\sigma>1$,
\begin{equation}\label{eq.resD}
  \|w_\sigma^{-\frac12}R_{|D|}(z)f\|_{L^2}
  \leq
  C\|w_\sigma^{\frac12} f\|_{L^2},
\end{equation}
uniformly in $z\in\C$.

We are ready to prove a corresponding estimate for the
resolvent of the perturbed operator
$$R(z)=(\sqrt{-\Delta+W}-z)^{-1},\qquad
W=a(x)\cdot\nabla+b(x),$$
following the same approach as in the preceding cases.

\begin{lemma}\label{lem.resDW}
  Consider the operator
  $-\Delta+W(x,D)\equiv-\Delta+a(x)\cdot \nabla+b_{1}(x)+b_{2}(x)$
  under the assumptions: the operator is selfadjoint,
  $b_{2}$ is real valued and nonnegative, and for some
  $\delta,\epsilon>0$ small enough, $\sigma>1$,
  \begin{equation}\label{eq.ipinversonde}
    \|\tau_\epsilon w_\sigma^{\frac12}
    a\|_{L^\infty}+
    \|\tau_\epsilon^2b_{1}\|_{L^\infty}< \delta,\qquad
    \|\tau_\epsilon^2b_{2}\|_{L^\infty}<\infty.
  \end{equation}
  Moreover assume that 0 is not a resonance for
  $-\Delta+b_{2}(x)$.
  Then the resolvent operator
  $R(z)=(\sqrt{-\Delta+W}-z)^{-1}$ satisfies
  \begin{equation}\label{eq.resDWstima}
    \|\tau_\epsilon^{-1}R(z)f\|_{L^2}
    \leq
    C\|\tau_\epsilon f\|_{L^2}.
  \end{equation}
  As a consequence, the perturbed wave
  flow $e^{it\sqrt{\Delta+W}}$ satisfies
  the smoothing estimate
  \begin{equation}\label{eq.katoonde}
  \|\tau_\epsilon^{-1}
  e^{it\sqrt{-\Delta+W}}f\|_{L^2L^2}\leq
  C\|f\|_{L^2}.
\end{equation}
\end{lemma}

\begin{proof}
We write for brevity
\begin{equation*}
  |D_{W}|=\sqrt{-\Delta+W(x,D)}.
\end{equation*}
By the (Phragm\'en-Lindel\"of) maximum principle, it is
sufficient to prove estimate \eqref{eq.resDWstima}
for real $z=\lambda$.
We notice that by the same arguments used
in the proof of Lemma \ref{lem.selfadj}, we have
\begin{equation*}
    \||D_{W}|g\|_{L^{2}}\simeq\|g\|_{\dot H^{1}};
\end{equation*}
thus for $\lambda\le0$ we can write
\begin{equation*}
    \|(|D_{W}|-\lambda)g\|^{2}_{L^{2}}
    =\||D_{W}|g\|^{2}_{L^{2}}+\lambda^{2}\|g\|^{2}_{L^{2}}-
        2\lambda(|D_{W}|g,g)_{L^{2}}
    \gtrsim\|g\|_{\dot H^{1}}
\end{equation*}
by the nonnegativity of $|D_{W}|$.
This implies for all $\lambda\le0$
\begin{equation*}
    \|R(\lambda)g\|_{\dot H^{1}}\lesssim\|g\|_{L^{2}},
\end{equation*}
whence by duality we have also
\begin{equation*}
    \|R(\lambda)g\|_{ L^{2}}\lesssim\|g\|_{\dot H^{-1}},
\end{equation*}
and interpolating we obtain
\begin{equation*}
    \|R(\lambda)g\|_{\dot H^{1/2}}\lesssim\|g\|_{\dot H^{-1/2}},
    \qquad\lambda\le0.
\end{equation*}
Now, using the Hardy's inequalities
\begin{equation*}
  \||x|^{-1/2}f\|_{L^{2}}\lesssim\|f\|_{\dot H^{1/2}}
  \qquad\text{or equivalently}\qquad
  \|f\|_{\dot H^{-1/2}}\lesssim\||x|^{1/2}f\|_{L^{2}}
\end{equation*}
we obtain the estimate
\begin{equation}\label{eq.AAAA}
    \||x|^{-1/2}R(\lambda)g\|_{L^2}\lesssim
    \||x|^{1/2}g\|_{L^{2}},\qquad\lambda\le0
\end{equation}
which implies \eqref{eq.resDWstima} for $z=-\lambda\le0$
(and is actually stronger).

Consider now $R(\lambda)$, $\lambda\ge0$; we use the
identity
\begin{equation*}
    R(\lambda)=(|D_{W}|-\lambda)^{-1}=
     2\lambda R_{W}(\lambda^{2})+(|D_{W}|+\lambda)^{-1}
\end{equation*}
where $R_W(\lambda)=(-\Delta+W-\lambda)^{-1}$.
The second term at the right hand side has already been
estimated, while the first one can be estimated using
\eqref{eq.estRW}, and this concludes the proof of
\eqref{eq.resDWstima}. The last inequality
\eqref{eq.katoonde} is an application of Kato's theorem
as usual.
\end{proof}

We conclude this section with a study of the operator
$\sqrt{-\Delta+1+W}$ associated with the perturbed
Klein-Gordon flow $e^{it\sqrt{-\Delta+1+W}}$.
In the free case $W=0$ the operator reduces to
$\DD=(1-\Delta)^{1/2}$ and its resolvent
\begin{equation*}
  R_{\DD}(z)=(\DD-z)^{-1}
\end{equation*}
can be handled in a similar way as $R_{|D|}$.

We start from estimates \eqref{eq.brv2} and \eqref{eq.tau}
which imply
\begin{equation*}
  {\zz}^{1/2}\|\tau_\epsilon^{-1}R_{0}(z)\|_{L^{2}}+
  \|w_{\sigma}^{-1/2}\nabla R_{0}(z)\|_{L^{2}}\lesssim
  \|\tau_{\epsilon}f\|_{L^{2}}.
\end{equation*}
As above, using the fact that $w_{\sigma}^{-1}$
is an $A_{2}$ weight, we can replace $\nabla$
with $|D|$ in the left hand side and hence
(recalling that $w_{\sigma}\lesssim \tau_{\epsilon}$)
we arrive at
\begin{equation*}
  {\zz}^{1/2}\|\tau_\epsilon^{-1}R_{0}(z)\|_{L^{2}}+
  \|w_{\sigma}^{-1/2}\DD R_{0}(z)\|_{L^{2}}\lesssim
  \|\tau_{\epsilon}f\|_{L^{2}}
\end{equation*}
which implies
\begin{equation}\label{eq:DD1}
  {\zz}^{1/2}\|\tau_\epsilon^{-1}R_{0}(z)\|_{L^{2}}+
  \|\tau_{\epsilon}^{-1}\DD R_{0}(z)\|_{L^{2}}\lesssim
  \|\tau_{\epsilon}f\|_{L^{2}}.
\end{equation}
Then using the identity
\begin{equation*}
  R_{\DD}(z)=(\DD+z)\cdot R_{0}(1-z^{2})
\end{equation*}
we obtain from \eqref{eq:DD1} the estimate
\begin{equation}\label{eq:DD2}
  \|\tau_{\epsilon}^{-1}R_{\DD}(z)f\|_{L^{2}}\lesssim
  \|\tau_{\epsilon}f\|_{L^{2}}.
\end{equation}
For the perturbed operator we have:

\begin{lemma}\label{lem.resDW2}
  Consider the operator
  $-\Delta+W(x,D)\equiv-\Delta+a(x)\cdot \nabla+b_{1}(x)+b_{2}(x)$
  under the assumptions: the operator is selfadjoint,
  $b_{2}$ is real valued and nonnegative, and for some
  $\delta,\epsilon>0$ small enough, $\sigma>1$,
  \begin{equation}\label{eq.ipinverskg}
    \|\tau_\epsilon w_\sigma^{\frac12}
    a\|_{L^\infty}+
    \|\tau_\epsilon^2b_{1}\|_{L^\infty}< \delta,\qquad
    \|\tau_\epsilon^2b_{2}\|_{L^\infty}<\infty.
  \end{equation}
  Moreover assume that 0 is not a resonance for
  $-\Delta+b_{2}(x)$.
  Then the resolvent operator
  $R(z)=(\sqrt{1-\Delta+W}-z)^{-1}$ satisfies
  \begin{equation}\label{eq.resDWkg}
    \|\tau_\epsilon^{-1}R(z)f\|_{L^2}
    \leq
    C\|\tau_\epsilon f\|_{L^2}.
  \end{equation}
  As a consequence, the perturbed Klein-Gordon
  flow $e^{it\sqrt{\Delta+1+W}}$ satisfies
  the smoothing estimate
  \begin{equation}\label{eq.katokg}
  \|\tau_\epsilon^{-1}
  e^{it\sqrt{-\Delta+1+W}}f\|_{L^2L^2}\leq
  C\|f\|_{L^2}.
\end{equation}
\end{lemma}

\begin{proof}
  Writing
  \begin{equation*}
    |D_{W}|=\sqrt{-\Delta+W(x,D)},\qquad
    \DDW=\sqrt{1-\Delta+W(x,D)}
  \end{equation*}
  we notice that
  \begin{equation*}
    \|\DDW f\|_{L^{2}}\simeq \|f\|_{L^{2}}+\||D_{W}|f\|_{L^{2}}
    \simeq\|f\|_{H^{1}}
  \end{equation*}
  by the same arguments used in the
  proof of Lemma \ref{lem.selfadj} and the identity
  \begin{equation*}
    \|\DDW f\|^{2}_{L^{2}}=
    ((1-\Delta+W)f,f).
  \end{equation*}
  Proceeding as in the proof of Lemma \ref{lem.resDW},
  we arrive at
  \begin{equation*}
      \|R(\lambda)g\|_{H^{1/2}}\lesssim\|g\|_{H^{-1/2}},
      \qquad\lambda\le0
  \end{equation*}
  for the resolvent $R=(\DDW-z)^{-1}$, and by Hardy
  inequality as before we obtain half of\eqref{eq.katokg}.

  For positive $\lambda$ we write
  \begin{equation*}
      R(\lambda)=(\DDW-\lambda)^{-1}=
       2\lambda R_{W}(\lambda^{2}-1)+(\DDW+\lambda)^{-1}
  \end{equation*}
  where $R_{W}(z)=(-\Delta+W-z)^{-1}$,
  and by \eqref{eq.estRW} and the first part of the proof
  we obtain \eqref{eq.resDWkg}. Kato's theorem gives
  \eqref{eq.katokg} as usual.
\end{proof}

\subsection{The magnetic Dirac operators}
\label{ssec.estresdir}

We now consider the resolvent of a perturbed Dirac operator
$\D+V(x)$. The proofs here will be short since we
shall rely on a few results proved in \cite{DAnconaFanelli06-c};
in particular, we recall
that if $V=V^{*}$ has a sufficienlty small
$L^{3,\infty}$ norm, hence under the
assumptions of Theorem \ref{thm.dirac}, the operator $\D+V$ is
self-adjoint on $L^{2}(\R^3,\C^4)$, with form domain
$H^1(\R^3,\C^4)$ and spectrum $\R$.
The same holds for the operator with nonzero
mass $\mathcal D+\beta+V$, but the
spectrum is $\mathbb{R}\setminus]-1,1[$.

Let us consider the massless case first.
We shall use the notations
\begin{equation}\label{eq.resdir}
  R_\D(z)=(-\D-zI_4)^{-1},\qquad
  R(z)=(-\D+V-zI_4)^{-1}
\end{equation}
where $I_4$ denotes the identity $4\times4$-matrix.
The following result is contained in
Proposition 3.6 of \cite{DAnconaFanelli06-c},
apart from the smoothing estimate which is a standard
consequence of Kato's theorem as above:


\begin{proposition}\label{prop.RV}
  Assume that the $4\times4$ matrix $V(x)=V^{*}(x)$ satisfies
  \begin{equation}\label{eq.ipVlap}
    \|w_{\sigma}V\|_{L^\infty}
    <\delta,
  \end{equation}
  for some $\delta$ sufficiently small and some $\sigma>0$.
  Then $\D+V$ satisfies the limiting absorption principle,
  i.e., the limit operators $R(\lambda\pm i0)$ exist
  in the topology of bounded operators from
  $L^{2}(w_{\sigma}^{1/2}dx)$ to
  $H^{1}(w_{\sigma}^{1/2}dx)$.
  Moreover the resolvent operator
  $R=(-\D+V-zI_4)^{-1}$ satisfies the estimate
  \begin{equation}\label{eq.estRV}
    \|w_{\sigma}^{-1/2}R(z)f\|_{L^2}\leq
    C\|w_{\sigma}^{1/2} f\|_{L^2},\qquad z\in\C.
  \end{equation}
  As a consequence, the Dirac flow satisfies the
  smoothing estimate
  \begin{equation}\label{eq.katodir3}
  \|w_{\sigma}^{-1/2}e^{it(\D+V)}f\|_{L^2L^2}
  \leq
  C\|f\|_{L^2}.
\end{equation}
\end{proposition}

We consider now the operators with mass
$\mathcal D+\beta$ and $\mathcal D+\beta+V$.
We shall use the notations
\begin{equation*}
  R_{\beta}(z)=(\mathcal D+\beta-zI_4)^{-1},
  \qquad
  R(z)=(\mathcal D+\beta+V-zI_4)^{-1}.
\end{equation*}
From the identities
\begin{equation*}
  \D^2=-\Delta I_4,\qquad
  (\mathcal D+\beta)^2=(1-\Delta)I_4,
\end{equation*}
we obtain the following representations in terms of
$R_{0}(z)=(-\Delta-z)^{-1}$
\begin{equation*}
  R_\D(z)=R_0(z^2)(\D+zI_4),\qquad
  R_{\beta}(z)=R_0(z^2-1)(\mathcal D+\beta+zI_4);
\end{equation*}
and hence we can write
\begin{equation}\label{eq.rapresdir}
  R_{\beta}(z)=R_0(z^2-1)\mathcal D+
  R_0(z^2-1)(\beta+zI_4).
\end{equation}
Then a straightforward application of estimate
\eqref{eq:DD1} gives
\begin{equation}\label{eq.estRH}
  \|\tau_{\epsilon}^{-1}R_{\beta}(z)f\|_{L^2}\leq C
  \|\tau_{\epsilon}f\|_{L^2}.
\end{equation}
uniformly in $z\in\C$.

In the perturbed case we can prove

\begin{proposition}\label{prop.RQ}
  Assume that the $4\times4$ matrix $V(x)=V^{*}(x)$  satisfies
  \begin{equation}\label{eq.ipQlap}
    \|\tau_\epsilon^2V\|_{L^\infty}<\delta,
  \end{equation}
  for some $\delta$ sufficiently small and $\epsilon>0$.
  Then the perturbed resolvent operator
  $R(z)=(\mathcal D+\beta+V-zI_4)^{-1}$ satisfies
  \begin{equation}\label{eq.estRQ}
    \|\tau_\epsilon^{-1}R(z)f\|_{L^2}\leq
    C\|\tau_\epsilon f\|_{L^2}.
  \end{equation}
  As a consequence, the flow $e^{it(\D+\beta+V)}$ satisfies the
  smoothing estimate
  \begin{equation}\label{eq.katodirmass}
  \|\tau_\epsilon^{-1}e^{it(\D+\beta+V)}f\|_{L^2L^2}
  \leq
  C\|f\|_{L^2}.
\end{equation}
\end{proposition}

\begin{proof}
The operator $V R_{\beta}(z)$ is bounded on
$L^{2}(\tau_{\epsilon}^{2}dx)$ with norm bounded by
$C \delta$ since
\begin{equation*}
  \|\tau_{\epsilon}VR_{\beta}(z)\|_{L^{2}}\leq
  \|\tau_{\epsilon}^{2}V\|_{L^{\infty}}
  \|\tau_{\epsilon}^{-1}R_{\beta}(z)\|_{L^{2}}\leq
  C \delta\|\tau_{\epsilon}f\|_{L^{2}}
\end{equation*}
by \eqref{eq.ipQlap} and \eqref{eq.estRH}.
Thus for $\delta$ small a Neumann expansion shows
that $(I+VR_{\beta}(z))^{-1}$ is well defined and
uniformly bounded on $L^{2}(\tau_{\epsilon}^{2}dx)$.
Hence the usual representation
\begin{equation*}
  R(z)=R_{\beta}(z)(I+VR_{\beta}(z))^{-1}
\end{equation*}
together with \eqref{eq.estRH} gives \eqref{eq.estRQ},
and \eqref{eq.katodirmass} follows.
\end{proof}

\section{Proof of the Strichartz Estimates}\label{sec.stri}

The method we shall follow is inspired by
\cite{RodnianskiSchlag04-a},
\cite{BurqPlanchonStalkerTahvildar03} and consists
in mixing Strich\-artz and smoothing estimates
for the free operator with smoothing estimates
for the perturbed operator. The main tool will
be the well-known Christ-Kiselev lemma \cite{ChKi}, which
can be stated as follows: given two Banach spaces $X,Y$ and
a bounded integral operator
$Tf=\int_{\mathbb{R}}K(t,s)f(s)ds$
from $L^{p}(\mathbb{R},X)$ to $L^{\tilde p}(\mathbb{R},Y)$,
then its truncated version
$Sf=\int_{0}^{t}K(t,s)f(s)ds$ is also
bounded on the same spaces, provided $p<\tilde p$
(the Hilbert transform being a trivial
counterexample for $p=\tilde p$). Thus to prove
an estimate of the form
\begin{equation*}
  \left\|\int_{0}^{t}e^{i(t-s)A}F(s)ds\right\|_{L^{p}_{t}L^{q}_{x}}
  \lesssim
  \|F\|_{L^{\tilde p}_{t}L^{\tilde q}_{x}}
\end{equation*}
it is sufficient to prove the untruncated estimate
\begin{equation*}
  \left\|\int_{\mathbb{R}}e^{i(t-s)A}F(s)ds\right\|_{L^{p}_{t}L^{q}_{x}}
  \lesssim
  \|F\|_{L^{\tilde p}_{t}L^{\tilde q}_{x}}
\end{equation*}
but only if $p<\tilde p$, which in particular excludes
endpoint-endpoint estimates where $p=\tilde p=2$.

\subsection{Schr\"odinger equation: proof of
Theorem \ref{thm.schrogrande}}
\label{ssec.strischro}

Notice that $u(t,x)=e^{it(-\Delta+W)}f$ satisfies the equation
$iu_{t}-\Delta u=-Wu$, hence we can write
\begin{equation*}
   e^{it(\Delta-W)}f=e^{it\Delta}f-
      \int_{0}^{t}e^{i(t-s)\Delta}W(x,D)u\, ds=
      I-II-III
\end{equation*}
with
\begin{equation*}
  I=e^{it\Delta}f,\qquad
  II=\int_{0}^{t}e^{i(t-s)\Delta}b(x) u\, ds, \qquad
  III=\int_{0}^{t}e^{i(t-s)\Delta}a(x)\cdot \nabla u\, ds.
\end{equation*}
The first term $I$ can be estimated directly
with standard Strichartz estimates:
\begin{equation}\label{eq.strilib}
  \|e^{it\Delta}f\|_{L^p_{t}L^q_{x}}\leq C\|f\|_{L^2}
\end{equation}
for any admissible couple $(p,q)$.
In order to estimate the second term we
resort to the Christ-Kiselev lemma and we
are reduced to estimate the untruncated integral
\begin{equation*}
  II_{1}=e^{it \Delta}\int e^{-is\Delta}b(x) u\, ds.
\end{equation*}
To this end we apply first the Strichartz estimates for the
free group, then the dual of the smoothing estimate from Proposition
\ref{prop.estreslap} in the special case $W=0$, i.e.,
\begin{equation*}
 \left\|\int e^{-is \Delta}F(s)ds\right\|_{L^{2}}\lesssim
 \|\tau_{\epsilon}F\|_{L^{2}L^{2}}
\end{equation*}
obtaining
\begin{equation*}
  \|II_{1}\|_{L^{p}L^{q}}\lesssim
  \left\|\int e^{-is \Delta}bu\,ds\right\|_{L^{2}}\lesssim
    \|\tau_{\epsilon}bu\|_{L^{2}L^{2}} \le
        \|\tau_{\epsilon}^{2}b\|_{L^{\infty}}
           \|\tau_{\epsilon}u\|_{L^{2}L^{2}}.
\end{equation*}
Then by assumption \eqref{eq.ab} and again the smoothing
estimate \eqref{eq.katoschro1} we conclude
\begin{equation}\label{eq:II}
  \|II\|_{L^{p}L^{q}}\lesssim \|f\|_{L^{2}}
\end{equation}
for any non-endpoint admissible couple $(p,q)$.

The last term $III$ is more delicate. We reduce it
as above to the untruncated form
\begin{equation*}
  III_{1}=e^{it \Delta}\int e^{-is \Delta}a \cdot \nabla u\,ds
\end{equation*}
and we apply to it the free Strichartz estimate and then the
following dual smoothing estimate:
\begin{equation}\label{eq:duals}
  \left\|\int e^{-is \Delta}F(s)ds\right\|_{L^{2}}\lesssim
     \| |D|^{-1/2}\chi F\|_{L^{2}L^{2}},
\end{equation}
valid for any function $\chi(x)\gtrsim w_{\sigma}(x)^{1/2}$.
Estimate \eqref{eq:duals} is proved
as follows: from \eqref{eq.brv2} we deduce, using the
fact that $w_{\sigma}$ is an $A_{2}$ weight, the equivalent
property
\begin{equation*}
  \|w_{\sigma}^{-1/2}|D|^{1/2}R_{0}(z)f\|_{L^{2}}\le
    C\|w_{\sigma}^{1/2}|D|^{-1/2}f\|_{L^{2}}
\end{equation*}
which implies, via Kato smoothing,
\begin{equation*}
  \|w_{\sigma}^{-1/2}|D|^{1/2}e^{it \Delta}f\|_{L^{2}L^{2}}\le
  \|f\|_{L^{2}}.
\end{equation*}
Since $\chi\gtrsim w_{\sigma}^{1/2}$ this gives also
\begin{equation*}
  \|\chi^{-1}|D|^{1/2}e^{it \Delta}f\|_{L^{2}L^{2}}\le
  \|f\|_{L^{2}}
\end{equation*}
and by duality we get \eqref{eq:duals}. Thus we arrive at
\begin{equation}\label{eq.partt}
  \|III_{1}\|_{L^{p}L^{q}}\lesssim
  \||D|^{-1/2}\chi a(x)\cdot \nabla u\|_{L^{2}L^{2}}
\end{equation}
Now assume we can prove the inequality
\begin{equation}\label{eq:commin}
  \||D|^{-1/2}\chi a(x)\cdot \nabla g\|_{L^{2}}\lesssim
  \|\tau_{\epsilon}^{-1}\nabla|D|^{-1/2}g\|_{L^{2}};
\end{equation}
then from \eqref{eq.partt} and the smoothing estimate
\eqref{eq.katoschro1b} we finally obtain
\begin{equation}\label{eq:III}
  \|III_{1}\|_{L^{p}L^{q}}\lesssim
  \|\tau_{\epsilon}^{-1}\nabla|D|^{-1/2}u\|_{L^{2}}\lesssim
  \|f\|_{L^{2}}
\end{equation}
which, together with \eqref{eq.strilib} and \eqref{eq:II},
concludes the proof of the Theorem.

It remains to check inequality \eqref{eq:commin}. We
rewrite it in the equivalent form
\begin{equation*}
  \||D|^{-1/2}\chi a(x)|D|^{1/2}
       \tau_{\epsilon}h\|_{L^{2}}\lesssim
  \|h\|_{L^{2}},
\end{equation*}
i.e., we need to prove that the operator
\begin{equation}\label{eq.opl2}
  T=|D|^{-1/2}\chi a(x)|D|^{1/2}\tau_{\epsilon}
\end{equation}
is bounded on $L^{2}$. We shall use the following lemma,
where we shall make use of several properties
of Lorentz spaces $L^{p,q}$ (see \cite{oneil}).

\begin{lemma}\label{lem:commut}
  Let $\alpha(x),\beta(x)$ be measurable functions on
  $\mathbb{R}^{n}$ such that for some $0<\delta<1/2$,
  some $\rho\in[0,n/2-\delta[$,
  and a radial function $\gamma(|x|)$, with
  $\gamma(s)$ decreasing, we have

  (i) $|\alpha(x)-\alpha(y)|\lesssim
       |x-y|^{1/2+\delta}(\gamma(|x|)+\gamma(|y|))$ and
    $\gamma\in L^{\frac{2n}{1+2 \rho+2 \delta},\infty}$

  (ii) $\alpha\beta\in L^{\infty}$,
    $|x|^{-\rho} \beta(x)\in L^{\infty}$
    and $|x|^{\rho}\gamma(|x|)\in L^{\frac{2n}{1+2 \delta},\infty}$

  \noindent
  Then the operator
  $T=|D|^{-1/2}\alpha(x)|D|^{1/2}\beta(x)$
  is bounded on $L^{2}$.

  The same result holds in the range $\rho\in[0,n/2+\delta[$
  if we replace (i) with

  (i') $|\alpha(x)-\alpha(y)|\le\langle x-y\rangle^{-2 \delta}
       |x-y|^{1/2+\delta}(\gamma(|x|)+\gamma(|y|))$ and
    $\gamma\in L^{\frac{2n}{1+2 \rho-2 \delta},\infty}$.
\end{lemma}

\begin{proof}

  Since $\alpha \beta$ is bounded, we can equivalently prove
  that the modified operator
  \begin{equation*}
    \widetilde{T}=T-\alpha \beta=
      |D|^{-1/2}\cdot[\alpha,|D|^{1/2}]\cdot\beta
  \end{equation*}
  is bounded on $L^{2}$. Moreover,
  by the Sobolev embedding in Lorentz spaces (proved e.g. by
  real interpolation)
  \begin{equation*}
    \||D|^{-1/2}g\|_{L^{2}}\lesssim
    \|g\|_{L^{\frac{2n}{n+1},2}}
  \end{equation*}
  it is sufficient to prove that the following
  reduced operator $S$ satisfies
  \begin{equation*}
    S=[\alpha,|D|^{1/2}]\cdot\beta:L^{\frac{2n}{n+1},2}\to L^{2}.
  \end{equation*}
  Now we observe that the commutator $[\alpha,|D|^{1/2}]$
  admits an explicit representation of the form
  \begin{equation*}
    [\alpha,|D|^{1/2}]f=c(n)
      \int_{\mathbb{R}^{n}}
      \frac{\alpha(x)-\alpha(y)}{|x-y|^{n+1/2}}\,f(y)dy
  \end{equation*}
  for a constant $c(n)$ depending only on the space dimension.
  Indeed, by standard Fourier transform techniques we see
  that
  \begin{equation*}
    [\alpha,|D|^{z}]f=c(z)
      \int_{\mathbb{R}^{n}}
      \frac{\alpha(x)-\alpha(y)}{|x-y|^{n+z}}\,f(y)dy
  \end{equation*}
  and this formula is valid for $\Re z<0$ under quite general
  assumptions on $\alpha$; moreover our assumptions
  show that the right hand side is a well defined and
  analytic function of $z$ for $\Re z<1/2+\delta$ (as proved
  below), hence by analytic continuation the representation is
  valid also in this larger region and in particular for $z=1/2$.

  In order to estimate $S$ we split it as $S=S_{1}+S_{2}$
  with
  \begin{equation*}
    S_{1}f=c\int_{|y|\ge2|x|}
    \frac{\alpha(x)-\alpha(y)}{|x-y|^{n+1/2}}\,\beta(y) f(y)dy
  \end{equation*}
  \begin{equation*}
    S_{2}f=c\int_{|y|\le2|x|}
    \frac{\alpha(x)-\alpha(y)}{|x-y|^{n+1/2}}\,\beta(y) f(y)dy
  \end{equation*}
  In the region $|y|\ge2|x|$ we deduce by assumption (i) that
  \begin{equation*}
    |\alpha(x)-\alpha(y)|\le 2|x-y|^{1/2+\delta} \gamma(|x|)
  \end{equation*}
  since $\gamma$ is decreasing; moreover we have
  $|x-y|\simeq |y|$, hence
  \begin{equation*}
    \left|
    \frac{\alpha(x)-\alpha(y)}{|x-y|^{n+1/2}}\,\beta(y) f(y)
    \right|\lesssim
    \gamma(|x|)
    \frac{|\beta(y)|}{|y|^{\rho}}
    \frac{|f(y)|}{|x-y|^{n-\rho-\delta}}\lesssim
    \gamma(|x|)
    \frac{|f(y)|}{|x-y|^{n-\rho-\delta}}
  \end{equation*}
  using (ii).
  Thus, by H\"older inequality in Lorentz spaces, we get
  \begin{equation*}
    \|S_{1}f\|_{L^{\frac{2n}{n+1},2}}\lesssim
    \|\gamma\|_{L^{\frac{2n}{1+2 \rho+2 \delta},\infty}}
    \left\|\int
    \frac{|f(y)|}{|x-y|^{n-\rho-\delta}}dy
    \right\|_{L^{\frac{2n}{n-2 \rho-2 \delta},2}}
  \end{equation*}
  (provided $\rho<n/2-\delta$)
  and by (i) and Young inequality we arrive at
  \begin{equation*}
    \|S_{1}f\|_{L^{\frac{2n}{n+1},2}}\lesssim
    \||y|^{-n+\rho+\delta}\|_{L^{\frac n {n-\rho-\delta},\infty}}
    \|f\|_{L^{2}}
  \end{equation*}
  which concludes.
  the estimate of the first piece $S_{1}$.

  In the region $|y|\le2|x|$, on the other hand, we can write
  \begin{equation*}
    \left|
    \frac{\alpha(x)-\alpha(y)}{|x-y|^{n+1/2}}\,\beta(y) f(y)
    \right|\lesssim
    \frac{|\beta(y)|}{|y|^{\rho}}
    \frac{|y|^{\rho}|\gamma(|y|/2)f(y)|}{|x-y|^{n-\delta}}
    \lesssim
    \frac{|y|^{\rho}|\gamma(|y|/2)f(y)|}{|x-y|^{n-\delta}}
  \end{equation*}
  so that by Young inequality
  \begin{equation*}
    \|S_{2}f\|_{L^{\frac{2n}{n+1},2}}\lesssim
    \left\|\int
    \frac{|y|^{\rho}|\gamma(|y|/2)f(y)|}{|x-y|^{n-\delta}}
    \right\|_{L^{\frac{2n}{n+1},2}}\lesssim
    \||y|^{\delta-n}\|_{L^{\frac{n}{n-\delta},\infty}}
    \|\gamma |y|^{\rho} f\|_{L^{\frac{2n}{n+1+2 \delta},2}}
  \end{equation*}
  and by H\"older inequality we get
  \begin{equation*}
    \|S_{2}f\|_{L^{\frac{2n}{n+1},2}}\lesssim
    \||y|^{\rho}\gamma\|_{L^{\frac{2n}{1+2 \delta},\infty}}
    \|f\|_{L^{2}}
  \end{equation*}
  and this concludes the proof under assumptions (i)-(ii).

  The case of assumptions (i')-(ii) is almost identical. No
  change is necessary in the estimate of $S_{2}f$, while
  for $S_{1}f$ it is sufficient to write
  \begin{equation*}
    \|S_{1}f\|_{L^{\frac{2n}{n+1},2}}\lesssim
    \|\gamma\|_{L^{\frac{2n}{1+2 \rho-2 \delta},\infty}}
    \left\|\int
    \frac{|f(y)|}{|x-y|^{n-\rho+\delta}}dy
    \right\|_{L^{\frac{2n}{n-2 \rho+2 \delta},2}}
  \end{equation*}
  which is true if $\rho<n/2+\delta$, and then proceed
  as above.
\end{proof}

Notice that if we restrict to the special choice
$\beta=|x|^{\rho}$, $\gamma(x)=\xx^{-\lambda}$,
$\alpha(x)=\chi(x)a(x)$, the following
conditions imply that (i), (ii), (i') are all satisfied:
\begin{equation}\label{eq:c1}
  0<\delta<\frac12,\qquad
  0\le\rho<\frac n2+\delta,\qquad
  \lambda\ge \frac12+\rho+\delta
\end{equation}
and
\begin{equation}\label{eq:c2}
  \xx^{\lambda}\chi(x)a(x)\in C^{1/2+\delta}
\end{equation}
(recall that $\|f\|_{C^{\mu}}=
\|f\|_{L^{\infty}}+\sup_{x\neq y}|x-y|^{-\mu}|f(x)-f(y)|$).
All conditions in (i), (ii), (i') are trivial to check
apart from H\"older continuity; actually we shall now
see that the following stronger inequality holds:
\begin{equation}\label{eq:c3}
  |\alpha(x)-\alpha(y)|\lesssim
       \min\{1,|x-y|\}^{1/2+\delta}(\xx^{-\lambda}+\yy^{-\lambda}).
\end{equation}
Indeed, when $|x-y|\ge1$ condition \eqref{eq:c3} follows from
$\xx^{\lambda}\chi(x)a(x)\in L^{\infty}$ which
is contained in \eqref{eq:c2}. When $|x-y|\le1$, we write
\begin{equation*}
  |\alpha(x)-\alpha(y)|\le A+B,
\end{equation*}
where
\begin{equation*}
    A=\chi(x)a(x)
       \xx^{\lambda}|\xx^{-\lambda}-\yy^{-\lambda}|,
\end{equation*}
and
\begin{equation*}
    B=\yy^{-\lambda}|\xx^{\lambda}\chi(x)a(x)-
        \yy^{\lambda}\chi(y)a(y)|.
\end{equation*}
Then we have directly from \eqref{eq:c2}
\begin{equation*}
  B \lesssim \yy^{-\lambda}|x-y|^{1/2+\delta}\le
      (\xx^{-\lambda}+\yy^{-\lambda})|x-y|^{1/2+\delta}
\end{equation*}
while for $A$ we use the elementary inequality
\begin{equation*}
  |\xx^{-\lambda}-\yy^{-\lambda}|\lesssim
    \sup_{\xi\in[x,y]}|\nabla\zz^{-\lambda}|_{z=\xi}\cdot|x-y|\lesssim
    (\xx^{-\lambda}+\yy^{-\lambda})|x-y|^{1/2+\delta}
\end{equation*}
together with the bound
$\xx^{\lambda}\chi(x)a(x)\in L^{\infty}$.

We can finally apply the lemma to the operator
\eqref{eq.opl2}; since
$\tau_{\epsilon}=|x|^{1/2-\epsilon}+|x|$ for $n\ge3$
and
$\tau_{\epsilon}=|x|^{1/2-\epsilon}+|x|^{1+\epsilon}$ for $n=2$,
by the above computation it is sufficient to
check conditions \eqref{eq:c1}, \eqref{eq:c2} for
$\rho=1/2-\epsilon$ and $\rho=1$
($\rho=1/2-\epsilon$ and $\rho=1+\epsilon$ in dimension 2).
We see that the choices $\delta=2 \epsilon$ and
$\lambda=1+3 \epsilon$ work in all cases, thus it is
sufficient to assume $\xx^{1+3 \epsilon}\chi(x)a(x)
\in C^{1/2+2 \epsilon}$ i.e. assumption
\eqref{eq.ipab2}. The proof is concluded.

\subsection{Wave and Klein-Gordon equations: proof of Theorems
\ref{thm.onde}, \ref{thm.kleingordon}}
\label{ssec.ondeklein}

Since $u(t,x)=e^{it\sqrt{-\Delta+W}}f$
solves the Cauchy problem
\begin{equation}\label{eq.ondemagn}
  \left\{\begin{array}{l}
      u_{tt}-\Delta u=-Wu
      \\
      u(0,x)=f(x)
      \\
      u_t(0,x)=i\left(\sqrt{-\Delta+W}
      \right)f(x),
    \end{array}\right.
\end{equation}
we have the alternative representation
\begin{equation}\label{eq.TT*ondemagn}
  e^{it\sqrt{-\Delta+W}}f=\cos(t|D|)f
  +i\frac{\sin(t|D|)}{|D|}\sqrt{-\Delta+W}f
  -\int_0^t\frac{\sin((t-s)|D|)}{|D|}Wuds.
\end{equation}
The first two terms satisfy the standard Strichartz estimates
for the free wave equation (see \eqref{eq.waveadmis}
in the Introduction, and recall also \eqref{eq.posit}).
For the third term we apply as usual the
Christ-Kiselev lemma and we are reduced to the untruncated
integral
\begin{equation*}
  \int\frac{\sin((t-s)|D|)}{|D|}Wuds=I+II
\end{equation*}
where, writing $c(x)=-\nabla \cdot a +b_{1}+b_{2}$,
\begin{equation*}
  I=\int\frac{\sin((t-s)|D|)}{|D|}\nabla\cdot(a(x) u)ds,\qquad
  II=\int\frac{\sin((t-s)|D|)}{|D|}c(x) uds.
\end{equation*}
Consider $I$; clearly, it is sufficient and actually stronger
to estimate the integral
\begin{equation*}
  I_{1}=|D|^{-1}e^{it |D|}\int e^{-is |D|} \nabla\cdot(a(x) u)ds=
     |D|^{-1}\nabla \cdot e^{it |D|}\int e^{-is |D|}a(x) uds.
\end{equation*}
To this end we recall the standard Strichartz estimate
\begin{equation}\label{eq.strm1}
    \|e^{it|D|}f\|_{L^p\dot H^{\frac1q-\frac1p-\frac12}_q}\lesssim
    \|f\|_{L^{2}}
\end{equation}
valid for any wave admissible couple $(p,q)$.
Moreover, the smoothing estimate \eqref{eq.katoonde}
holds also in the free case $W \equiv0$
\begin{equation}\label{eq.katoonde2}
  \|\tau_\epsilon^{-1}e^{it|D|}f\|_{L^2L^2}
  \lesssim\|f\|_{L^{2}}
\end{equation}
and by duality is equivalent to
\begin{equation}\label{eq.katoonde3}
  \left\|\int e^{-is|D|}F(s)ds\right\|_{L^{2}}
  \lesssim\|\tau_\epsilon F\|_{L^{2}L^{2}}.
\end{equation}
Applying \eqref{eq.strm1} and \eqref{eq.katoonde3} to $I_{1}$
we obtain, since the Riesz operators are bounded in all
$L^{p}$ with $1<p<\infty$,
\begin{equation*}
  \|I_{1}\|_{L^p\dot H^{\frac1q-\frac1p-\frac12}_q}\lesssim
  \|\tau_\epsilon a(x)u\|_{L^{2}L^{2}}\le
  \|\tau_\epsilon^{2} a(x)\|_{L^{2}}
    \|\tau_\epsilon^{-1} u\|_{L^{2}L^{2}}.
\end{equation*}
Using again the smoothing estimate \eqref{eq.katoonde}
and assumption \eqref{eq.ipabonde}, we conclude
\begin{equation*}
  \|I_{1}\|_{L^p\dot H^{\frac1q-\frac1p-\frac12}_q}\lesssim
    \|f\|_{L^{2}}.
\end{equation*}
Consider now the second term $II$, or more generally
\begin{equation*}
  II_{1}=e^{it|D|}\int |D|^{-1}e^{-is|D|}c(x)uds.
\end{equation*}
Proceeding as in \cite{BurqPlanchonStalker04-a}, we shall use the
following estimate from \cite{benar} (see also \cite{hosh})
\begin{equation*}
  \||x|^{-1}|D|^{-1}e^{it|D|}f\|_{L^{2}L^{2}}\lesssim
    \|f\|_{\dot H^{-1/2}}
\end{equation*}
in the dual form:
\begin{equation}\label{eq:hosh}
  \left\|\int |D|^{-1} e^{-is|D|} F(s)ds\right\|_{\dot H^{1/2}}
  \lesssim
  \||x|F\|_{L^{2}L^{2}}.
\end{equation}
Then, applying the Strichartz estimate for the wave
equation \eqref{eq.strm1} in the form
\begin{equation*}
  \|e^{it|D|}f\|_{L^p\dot H^{\frac1q-\frac1p}_q}\lesssim
  \|f\|_{\dot H^{\frac12}}
\end{equation*}
followed by \eqref{eq:hosh}, we obtain
\begin{equation*}
  \|II_{1}\|_{L^p\dot H^{\frac1q-\frac1p-\frac12}_q}\lesssim
  \||x|c(x)u\|_{L^{2}L^{2}}\lesssim
  \||x|\tau_{\epsilon}c(x)\|_{L^{\infty}}
    \|\tau_{\epsilon}^{-1}u\|_{L^{2}L^{2}}.
\end{equation*}
Recalling assumption \eqref{eq.ipabonde}
and the smoothing estimate
\eqref{eq.katoonde} we finally obtain
\begin{equation*}
  \|II_{1}\|_{L^p\dot H^{\frac1q-\frac1p-\frac12}_q}\lesssim
  \|f\|_{L^{2}}
\end{equation*}
which concludes the proof of Theorem \ref{thm.onde}.

The proof of Theorem \ref{thm.kleingordon} is completely
analogous, using the Strichartz estimate for the free
equation
\begin{equation*}
  \|e^{it\DD}f\|_{L^p H^{\frac1q-\frac1p-\frac12}_q}\lesssim
  \|f\|_{L^{2}},
\end{equation*}
which is valid for all Schr\"odinger admissible couple
$(p,q)$, and the following estimate from \cite{benar}:
\begin{equation*}
  \|\xx^{-1}e^{it\DD}f\|_{L^{2}L^{2}}\lesssim\|f\|_{L^{2}}
\end{equation*}
which implies by duality
\begin{equation*}
  \left\|\int e^{-is\DD} F(s)ds\right\|_{L^{2}}
  \lesssim
  \|\xx F\|_{L^{2}L^{2}}
\end{equation*}
and hence also
\begin{equation*}
  \left\|\int\DD^{-1} e^{-is\DD} F(s)ds\right\|_{H^{1/2}}
  \lesssim
  \|\xx F\|_{L^{2}L^{2}}.
\end{equation*}
This estimate replaces \eqref{eq:hosh}
in the above computation.

\subsection{Dirac equation: proof of Theorems \ref{thm.dirac},
\ref{thm.diracmass}}
\label{ssec.stridir}

As proved in the Appendix, the
Strichartz estimate for the free massless Dirac equation
is the following:
\begin{equation}\label{eq.stridirbis}
  \|e^{it\mathcal D}f\|_{L^p\dot H^{\frac1q-\frac1p-\frac12}_q}
  \lesssim \|f\|_{L^2}
\end{equation}
for any wave admissible couple  $(p,q)$.
On the other hand, as a special case of the smoothing
estimate \eqref{eq.katodir3}, we have
\begin{equation}\label{eq.katodir}
  \|w_\sigma^{-\frac12}e^{it\mathcal D}f\|_{L^2L^2}
  \lesssim\|f\|_{L^2}
\end{equation}
and by duality we obtain
\begin{equation}\label{eq.katodir2}
  \left\|\int e^{-is\mathcal D}
  F(s)ds\right\|_{L^2}
  \lesssim\|w_\sigma^{\frac12}F\|_{L^2L^2}.
\end{equation}

Consider now the perturbed Dirac flow
$u=e^{it(\mathcal D+V)}f$.
An alternative representation of $u$ is the following:
\begin{equation}\label{eq.rapprdir}
  u(t,x)=e^{it\mathcal D}f-e^{it\mathcal{D}}
  \int_0^te^{-is\mathcal D}Vu(s)ds.
\end{equation}
The term $e^{it\mathcal D}f$ satisfies the free
Strichartz estimates \eqref{eq.stridir}; in order
to estimate the Duhamel term as usual we apply the
Christ-Kiselev lemma and switch to the untruncated
integral. Then, using
\eqref{eq.stridir2}, \eqref{eq.katodir2} and H\"older
inequality, we have
\begin{eqnarray}
  \left\|e^{it\mathcal{D}}
     \int e^{-is\mathcal D}Vuds
  \right\|_{L^p\dot H^{\frac1q-\frac1p-\frac12}_q}
  & \lesssim &
  \left\|\int e^{-is\mathcal D}
  Vu ds\right\|_{L^2} \nonumber
  \\
  & \lesssim &
  \|w_\sigma^{\frac12} Vu\|_{L^2L^2}
  \leq
  \|w_\sigma V\|_{L^\infty}
  \cdot\|w_\sigma^{-\frac12}u\|_{L^2L^2}.
  \label{eq.integral}
\end{eqnarray}
Recalling the smoothing estimate \eqref{eq.katodir3}
we obtain
\begin{equation*}
  \left\|e^{it\mathcal{D}}
     \int e^{-is\mathcal D}Vuds
  \right\|_{L^p\dot H^{\frac1q-\frac1p-\frac12}_q}
  \lesssim
  \|f\|_{L^{2}}
\end{equation*}
and this completes the proof of \ref{thm.dirac}.

The proof of Theorem \ref{thm.diracmass} is completely
analogous.

\appendix
\section{Strichartz estimates for the free flows}
\label{sec.appendix}

Strichartz estimates for the free Schr\"odinger and wave
equations are well known, see the Introduction for the
precise statements.
It is less easy to find in the literature optimal results
for Klein-Gordon and Dirac equations. Hence we devote
this appendix to a quick proof of the estimates in these cases.

The massless Dirac flow is trivial since it can be
reduced to the wave equation:

\begin{proposition}\label{prop.stridir}
    Let $n=3$. The following Strichartz estimates hold:
    \begin{equation}\label{eq.stridir}
        \|e^{it\D}f\|_{L^p\dot H^{\frac1q-\frac1p-\frac12}_q}\lesssim
        \|f\|_{L^{2}}
    \end{equation}
    for any wave admissible couple $(p,q)$.
\end{proposition}

\begin{proof}
 By the identity
\begin{equation*}
  (i\partial_t+\mathcal D)(i\partial_t-\mathcal D)=
  -\square I_4,
\end{equation*}
we obtain that $u(t,x)=e^{it\D}f$ satisfies
the Cauchy problem
\begin{equation}\label{eq.wave2}
    \left\{\begin{array}{l}
    u_{tt}-\Delta I_4 u=0 \\
    u(0,x)=f(x) \\
    u_t(0,x)=i\mathcal Df(x)\end{array}\right.
\end{equation}
and hence each component of $u$ satisfies the same
Strichartz estimates as for the 3D wave equation.
\end{proof}

The Klein-Gordon and massive Dirac equations need some work.
We begin by the free Klein-Gordon flow $u=e^{it\DD}f$.
We shall apply a precise stationary phase
result due to H\"ormander \cite{Horm}:

\begin{lemma}\label{lem.horm}
    Assume that $\phi:\R^n\to\R$ has a Fourier transform
    $\widehat\phi\in\mathcal C^\infty$ with the
    decay property
    \begin{equation}\label{eq.ipphi}
        \left|D^\alpha\widehat\phi(\xi)\right|
        \leq C_\alpha\langle\xi\rangle^{-\frac n2-1-|\alpha|}
        \qquad \forall\xi\in\R^n,\ \alpha\in \mathbb{N}^{n}.
    \end{equation}
    Then the following estimate holds: for some $C>0$,
    \begin{equation}\label{eq.stima}
        \left|e^{it\langle D\rangle}\phi\right|
        \leq C(|t|+|x|)^{-\frac n2}.
    \end{equation}
\end{lemma}

Now, using an inhomogeneous dyadic decomposition
$\{\psi_{0},\varphi_j(D)\}_{j\geq1}$
with the usual properties:
$\psi_{0}(\xi)$ supported in $B(0,1)$,
$\varphi_{0}(\xi)=\psi_{0}(\xi/2)-\psi_{0}(\xi)$,
\begin{equation*}
    \varphi_j(\xi)=\varphi_0(2^{-j}\xi),\qquad
    \psi_{0}+\sum_{j\geq1}\varphi_j=1
\end{equation*}
we can localize the estimate as follows:

\begin{lemma}\label{lem.freq}
    The flow $e^{it\langle D\rangle}f$ satisfies the
    localized dispersive estimate
    \begin{equation}\label{eq.disploc}
        |e^{it\langle D\rangle}\varphi_j(D)f|\leq
        C|t|^{-\frac n2}2^{j(\frac n2+1)}
        \|\varphi_j(D)f\|_{L^1},
    \end{equation}
    for each $t\in\R$, $x\in\R^n$, $j\geq0$
    and some $C>0$; here
    $\widetilde\varphi_{j}$ denotes
    $\varphi_{j-1}+\varphi_{j}+\varphi_{j+1}$,
    with $\varphi_{-1}=0$.
\end{lemma}

\begin{proof}
We can write
\begin{equation*}
    e^{it\langle D\rangle}\varphi_j(D)f
    =e^{it\langle D\rangle}\langle
    D\rangle^{-\frac n2-1}\langle
    D\rangle^{\frac n2+1}\varphi_j(D)f=
    e^{it\langle D\rangle}\mathcal F^{-1}
    \left(\langle\xi\rangle^{-\frac n2-1}\right)*
    \left(\varphi_j(D)f\right),
\end{equation*}
where $\mathcal F^{-1}$ denotes the inverse Fourier transform.
Then, applying Lemma \ref{lem.horm} with
$\phi=\mathcal F^{-1}\left(\langle\xi
\rangle^{-\frac n2-1}\right)$, we obtain
\begin{equation}\label{eq.trucco2}
    \left|e^{it\langle D\rangle}\varphi_j(D)f
    \right|\leq C|t|^{-\frac n2}\|\langle
    D\rangle^{\frac n2+1}\varphi_j(D)f\|_{L^1}.
\end{equation}
Since
\begin{equation*}
    \langle D\rangle^{\frac n2+1}
    \varphi_j(D)f=\mathcal F^{-1}\left(
    \langle\xi\rangle^{\frac n2+1}\varphi_j(\xi)
    \right)*f,
\end{equation*}
Young inequality gives
\begin{equation}\label{eq.trucco3}
    \|\langle D\rangle^{\frac n2+1}
    \varphi_j(D)f\|_{L^1}\leq
    \|\mathcal F^{-1}\left(
    \langle\xi\rangle^{\frac n2+1}\varphi_j(\xi)
    \right)\|_{L^1}\|f\|_{L^1}.
\end{equation}
Notice that we can replace in this computation
$f$ with $\widetilde\varphi_{j}(D)f$ since
$\varphi_{j}(D)\widetilde\varphi_{j}(D)=\varphi_{j}(D)$.
Thus to conclude the proof it is sufficient to get
the following estimate:
\begin{equation}\label{eq.ultima}
    \|\mathcal F^{-1}\left(
    \langle\xi\rangle^{\frac n2+1}\varphi_j(\xi)
    \right)\|_{L^1}\leq C2^{j(\frac n2+1)}.
\end{equation}
Using the scaling operators
$S_{\lambda}\phi(x)=\phi(\lambda x)$,
we can write
\begin{eqnarray*}
    \mathcal F^{-1}\left(
    \langle\xi\rangle^{\frac n2+1}\varphi_j(\xi)
    \right) & = & \mathcal F^{-1}\left(
    \langle\xi\rangle^{\frac n2+1}S_{2^{-j}}
    \varphi_0(\xi)\right) \\ \  & = &
    2^{j(\frac n2+1)}2^{jn}S_{2^j}
    \mathcal F^{-1}\left(
    (2^{-2j}+|\xi|^2)^{\frac n2+1}\varphi_0(\xi)
    \right)
\end{eqnarray*}
and hence
\begin{equation*}
    \|\mathcal F^{-1}\left(
    \langle\xi\rangle^{\frac n2+1}\varphi_j(\xi)
    \right)\|_{L^1}\leq2^{j(\frac n2+1)}
    \|\mathcal F^{-1}\left(
    (2^{-2j}+|\xi|^2)^{\frac n2+1}\varphi_0(\xi)
    \right)\|_{L^1}.
\end{equation*}
Moreover, multiplying and dividing by $\xx^{2m}$
for some integer $m$, we obtain
\begin{eqnarray}
    \|\mathcal F^{-1}\left(
    (2^{-2j}+|\xi|^2)^{\frac n2+1}\varphi_0(\xi)
    \right)\|_{L^1} & \leq & C
    \|\xx^{2m}\mathcal F^{-1}\left(
    (2^{-2j}+|\xi|^2)^{\frac n2+1}\varphi_0(\xi)
    \right)\|_{L^2}\nonumber \\ \  & = & C
    \|(1-\Delta)^m\left(
    (2^{-2j}+|\xi|^2)^{\frac n2+1}\varphi_0(\xi)
    \right)\|_{L^2},\label{eq.ultima3}
\end{eqnarray}
provided
\begin{equation}\label{eq.m}
    m>\frac n4.
\end{equation}
We shall choose $m$ as the smallest integer satisfying
\eqref{eq.m}.
We are interested in the growth with respect to $j$
of the quantity
\begin{equation*}
  I:=(1-\Delta)^m\left(
    (2^{-2j}+|\xi|^2)^{\frac n2+1}\varphi_0(\xi)
    \right).
\end{equation*}
When $n$ is even,
$\left(2^{-2j}+|\xi|^2\right)^{\frac n2+1}$
is a polynomial, and hence we obtain
\begin{equation*}
  \|I\|_{L^2}\leq C\|\varphi_0\|_{L^2}
\end{equation*}
with $C$ independent of $j$. When $n$ is odd,
it is clear that almost all the terms in the
expansion of $I$ are uniformly bounded in $j$, apart
from the (possibly) worst one
\begin{equation*}
  II=\Delta^m\left(2^{-2j}+
  |\xi|^2\right)^{\frac n2+1}.
\end{equation*}
We have the two possibilities
\begin{equation*}
    n=4k+3\qquad\text{or}\qquad n=4k+1,
\end{equation*}
with $m=k+1$. If $n=4k+3$, we have
\begin{equation*}
  |II|\simeq\left|D^{2k+2}\left(
    (2^{-2j}+|\xi|^2)^{2k+\frac52}\right)
    \right|
\end{equation*}
which expands in a sum of bounded terms.
If $n=4k+1$, we have
\begin{equation*}
  |II| \simeq\left|D^{2k+2}\left(
  (2^{-2j}+|\xi|^2)^{2k+\frac32}\right)
  \right| \lesssim \left(2^{-2j}+|\xi|^2\right)^{-1/2}
  |\xi|^{2k+2}
  +\text{ bounded terms,}
\end{equation*}
and also in this case we have a uniform bound in $j$.
In conclusion, we have proved that
\begin{equation*}
  \|(1-\Delta)^m\left((2^{-2j}+|\xi|^2)^{\frac n2+1}
  \varphi_0(\xi)\right)\|_{L^2}\leq C,
\end{equation*}
for some $C>0$, which implies \eqref{eq.ultima}, and
the proof is complete.
\end{proof}

\begin{remark}
  By interpolation between estimate \eqref{eq.disploc}
  and the localized $L^2$ conservation
  \begin{equation}\label{eq.l2cons}
    \|e^{it\langle D\rangle}\varphi_j(D)f\|_{L^2}
    \leq\|\varphi_j(D)f\|_{L^2},
  \end{equation}
  we obtain the following $L^q-L^{q'}$ decay estimates:
  \begin{equation}\label{eq.displocp}
    \|e^{it\langle D\rangle}\varphi_j(D)f\|_{L^q}
    \leq C|t|^{-\frac n2+\frac nq}2^{j(\frac n2+1)
    (1-\frac2p)}\|\widetilde\varphi_j(D)f\|_{L^{q'}}
  \end{equation}
  for any $q\geq2$ with $1/q+1/q'=1$.
\end{remark}

Starting from estimates \eqref{eq.displocp} and
using the standard techniques
of \cite{GinibreVelo95-generstric},
\cite{KeelTao98-endpoinstric}, in particular
the abstract Theorem 10.1 of \cite{KeelTao98-endpoinstric},
we obtain the full set of estimates including the
endpoint case:

\begin{theorem}\label{thm.striKG}
  The Klein-Gordon flow $u=e^{it\DD}f$ satisfies
  the Strich\-artz estimates
  \begin{equation}\label{eq.striKG}
    \|e^{it\DD}f\|_{L^pH^{\frac1q-\frac1p-\frac12}_q}\lesssim
    \|f\|_{L^{2}}
  \end{equation}
  for any Schr\"odinger admissible couple
  $(p,q)$.
\end{theorem}

Finally, the Dirac equation with mass can be
handled in a similar way to Proposition \ref{prop.stridir}:

\begin{proposition}\label{prop.stridir2}
    Let $n=3$. The following Strichartz estimates hold:
    \begin{equation}\label{eq.stridirmass}
      \|e^{it(\DD+\beta)}f\|_{L^pH^{\frac1q-\frac1p-\frac12}_q}
      \lesssim
      \|f\|_{L^{2}},
    \end{equation}
    for any Schr\"odinger admissible couple
    $(p,q)$.
\end{proposition}

\begin{proof}
  As in the proof of Proposition \ref{prop.stridir},
  by the identity
\begin{equation*}
  (i\partial_t+(\mathcal D+\beta))(i\partial_t-(\mathcal D+\beta))
  =(-\square-1)I_4
\end{equation*}
we obtain that each component of $u$ solves a Klein-Gordon
equation with initial data $f$ and $(\mathcal D+\beta)f$. Thus
estimate \eqref{eq.stridirmass} follows immediately from the Strichartz
estimates for the Klein-Gordon equation in space dimension $n=3$.
\end{proof}

\end{document}